\documentclass[a4paper,reqno,11pt,draft]{amsart}
\usepackage{amsmath,amssymb,amsfonts,amscd}

\usepackage{verbatim}

\title{Bousfield
localization on formal schemes}

\author[L. Alonso]{Leovigildo Alonso Tarr\'{\i}o}
\address[L. A. T.]{Departamento de \'Alxebra\\
Facultade de Matem\'a\-ticas\\
Universidade de Santiago de Compostela\\
E-15782  Santiago de Compostela, SPAIN}
\email{leoalonso@usc.es}

\author[A. Jerem\'{\i}as]{Ana Jerem\'{\i}as L\'opez}
\address[A. J. L.]{Departamento de \'Alxebra\\
Facultade de Matem\'a\-ticas\\
Universidade de Santiago de Compostela\\
E-15782  Santiago de Compostela, SPAIN}
\email{jeremias@usc.es}

\author[M. J. Souto]{Mar\'{\i}a~Jos\'e Souto Salorio}
\address[M. J. S. S.]{Facultade de Inform\'atica, Campus de Elvi\~na\\
Universidade da Coru\~{n}a\\
E-15071  A Coru\~{n}a, SPAIN}
\email{mariaj@udc.es}

\thanks{L.A.T. and A.J.L. partially supported by Spain's MCyT  and E.U.'s
FEDER research project BFM2001-3241, supplemented by Xunta de Galicia grant
PGDIT 01PX120701PR}

\subjclass{14F99 (primary); 14F05, 18E30 (secondary)}

\date{December 7, 2003}

\theoremstyle{plain}
\newtheorem{thm}{Theorem}[section]
\newtheorem{lem}[thm]{Lemma}
\newtheorem{cor}[thm]{Corollary}
\newtheorem{prop}[thm]{Proposition}

\theoremstyle{remark}
\newtheorem*{rem}{Remark}

\theoremstyle{definition}

\newtheorem*{ex}{Example}
\newtheorem*{ack}{Acknowledgements}
\newtheorem{cosa}[thm]{}

\newcommand{\CA}{\mathcal{A}}
\newcommand{\CB}{{\mathcal B}}
\newcommand{\CC}{{\mathcal C}}

\newcommand{\CE}{{\mathcal E}}
\newcommand{\CF}{{\mathcal F}}
\newcommand{\CG}{{\mathcal G}}
\newcommand{\CH}{{\mathcal H}}
\newcommand{\CI}{{\mathcal I}}
\newcommand{\CJ}{{\mathcal J}}
\newcommand{\CK}{{\mathcal K}}
\newcommand{\CL}{\mathcal{L}}
\newcommand{\CM}{\mathcal{M}}
\newcommand{\CN}{\mathcal{N}}
\newcommand{\CO}{\mathcal{O}}
\newcommand{\CP}{\mathcal{P}}
\newcommand{\CQ}{\mathcal{Q}}

\newcommand{\CS}{\mathcal{S}}

\newcommand{\FU}{\mathfrak{U}}
\newcommand{\FX}{\mathfrak{X}}

\newcommand{\D}{\boldsymbol{\mathsf{D}}}
\newcommand{\K}{\boldsymbol{\mathsf{K}}}
\newcommand{\LL}{\boldsymbol{\mathsf{L}}}
\newcommand{\R}{\boldsymbol{\mathsf{R}}}
\newcommand{\T}{\boldsymbol{\mathsf{T}}}

\newcommand{\ts}{\mathsf{t}}

\newcommand{\pf}{\mathsf{cp}}
\newcommand{\qc}{\mathsf{qc}}
\newcommand{\qct}{\mathsf{qct}}

\newcommand{\PR}{\mathbf{P}}

\newcommand{\NN}{\mathbb{N}}
\newcommand{\ZZ}{\mathbb{Z}}

\newcommand{\ip}{{\mathfrak p}}

\newcommand{\kk}{\kappa}
\newcommand{\gmm}{\gamma\,}

\newcommand{\holim}[1]{\begin{array}[t]{c} {\rm holim}\\[-7.5 pt]
 {\longrightarrow} \\[-7.5 pt] {\scriptstyle {#1}} \end{array}}
\newcommand{\dirlim}[1]{\begin{array}[t]{c} {\rm lim}\\[-7.5 pt]
 {\longrightarrow} \\[-7.5 pt] {\scriptstyle {#1}} \end{array}}

\newcommand{\smdirlim}[1]{\begin{array}[t]{c} \scriptstyle {\rm lim}\\
[-8 pt] {\rightarrow} \\[-8 pt] {\scriptscriptstyle {#1}} \end{array}}

\newcommand{\lto}{\longrightarrow}
\newcommand{\xto}{\xrightarrow}
\newcommand{\ot}{\leftarrow}

\newcommand{\inc}{\hookrightarrow}
\newcommand{\iso}{\tilde{\to}}
\newcommand{\osi}{\tilde{\ot}}
\newcommand{\liso}{\tilde{\lto}}

\newcommand{\imp}{\Rightarrow}
\newcommand{\dimp}{\Leftrightarrow}
\newcommand{\st}{\mathsf{Ho}\mathbf{Sp}}

\DeclareMathOperator{\Hom}{Hom}
\DeclareMathOperator{\shom}{\CH\mathit{om}}
\DeclareMathOperator{\dhom}{\boldsymbol{\CH}\mathsf{om}}

\DeclareMathOperator{\spec}{Spec}
\DeclareMathOperator{\spf}{Spf}
\DeclareMathOperator{\supp}{Supp}
\DeclareMathOperator{\suph}{Supph}
\DeclareMathOperator{\id}{id}

\DeclareMathOperator{\loc}{\mathsf{Loc}}
\DeclareMathOperator{\Pa}{\mathsf{P}}
\DeclareMathOperator{\thik}{\mathsf{Th}}

\newcommand{\ie}{{\it i.e.} }

\begin{document}

\begin{abstract} Let $(\FX, \CO_\FX)$ be a noetherian formal scheme
and consider $\D_\qct(\FX)$ its derived category of sheaves with 
quasi-coherent torsion homology. We show that there is a bijection between
the set of rigid (\ie $\otimes$-ideals) localizing subcategories of
$\D_\qct(\FX)$ and subsets in $\FX$, generalizing previous work by Neeman. 
If moreover $\FX$ is separated, the associated localization and
acyclization functors are described in certain cases. When $Z \subset \FX$ is
a stable for specialization subset, its associated acyclization is
$\R\varGamma_Z$. When $X$ is an scheme, the corresponding localizing
subcategories are generated by perfect complexes and we recover Thomason's
classification of thick subcategories. On the other hand, if $Y \subset \FX$
is generically stable, we  show that the associated localization functor is
$\dhom^\cdot_{\FX}(\R{}\varGamma_{\FX \setminus Y} \CO'_\FX, \CG))$.
\end{abstract}

\maketitle

\section*{Introduction}

The techniques of localization have a long tradition in several areas of
mathematics. They have the virtue of concentrating our attention on some part
of the structure in sight allowing us to handle more manageable pieces
of information. One of the clear examples of this technique is the
localization in algebra where one studies a module centering the attention
around a point of the spectrum of the base ring, \ie a prime ideal. The
idea was transported to topology by Adams and later Bousfield proved that
there are plenty of localizations in stable homotopy. In the past decade it
became clear that one could successfully transpose homotopy techniques to the
study of derived categories (over rings and schemes). In particular, in
our previous work, we have shown that for the derived category
of a Grothendieck category we also have plenty of localizations. In that
paper, \cite{AJS}, we applied the result to the existence of unbounded
resolutions and we hinted that, in the case of the derived category of
quasi-coherent sheaves over a nice scheme, there should be a connection 
between localizations in the derived category and the geometric structure of
the underlying space.

The present paper realizes that goal extending the work of Neeman
\cite[Theorem 3.3]{Nct} who classified all Bousfield localizations in the
derived category of modules over a noetherian ring $\D(R)$ to the
classification of the Bousfield localizations of the derived category of
sheaves with quasi-coherent \emph{torsion} homology over a noetherian formal
scheme (Theorem \ref{ThclasB}). This category is a basic ingredient in
Grothendieck duality \cite{dfs}. Also, if the formal scheme is just an
\emph{usual} noetherian scheme it gives the derived category of sheaves with
quasi-coherent homology. Thus we obtain an analog of the \emph{chromatic
tower} in stable homotopy for these kind of schemes and formal schemes. It is
clear that the monoidal structure of the derived category is an  essential
part of the cohomological formalism. In fact, to get the classification we
were forced to consider only \emph{rigid} localizing subcategories. This
means, roughly speaking, that the localizing subcategory is an ideal in the
monoidal sense (see \S 3). This condition is needed in order to have
compatibility with open sets. It holds for all localizing subcategories in
the affine case, that is why it was not considered by Neeman. 

The classification theorem is more useful if the localization functor
associated to a subset of the formal scheme can be expressed in geometrically
meaningful terms with respect to this subset. This can be done for
noetherian separated formal schemes under certain conditions over the subset.
The most rich case is the case of stable for specialization subsets (that
recover the classical system of supports). They provide localizations that
have the property of being compatible with the tensor product. They are also
characterized by being associated to a right-derived functor and they
correspond to the smashing localizations of topologists. All of this is
contained in Theorem \ref{tenscomp}. These kind of localizations correspond to
Lipman's notion of idempotent pairs \cite{Lrtfs}. The associated localizing
subcategory is characterized in terms of homological support (Theorem
\ref{suph}). With this tool at hand we see that our classification of tensor
triangulated categories (or smashing localizations) agrees with Thomason
classification of thick $\otimes$-subcategories of the derived category
quasi-coherent sheaves \cite{th}, when both make sense \ie for a noetherian
separated scheme.

The dual notion of \emph{tensor compatible} is that of \emph{Hom compatible}
localizations. They correspond to stable for generalization subsets, which
are complementary of stable for specialization subsets. The Hom compatible
localizations can be described via a certain formal duality relation with the
tensor compatible localization associated to its complementary subsets
(Theorem \ref{homcomp}). If the stable for generalization subset is an open
set, the localization functor agrees with the left derived of a completion.
This relates the results of \cite{AJL1} to this circle of ideas.

While our work does not exhaust all the possible questions about these
topics we believe that it can be useful for the current program of extracting
information on a space looking at its derived category. 

Now, let us describe briefly the contents of the paper. The first section
recalls the concepts and notations used throughout and we give a detailed
overview of the symmetric closed structure in the derived categories we are
going to consider. In the next section, we specify the relationship between
cohomology with supports and the algebraic version defined in terms of ext
sheaves. We make a detailed study of the cohomology with respect to a system
of supports in the case of a formal scheme and interpret the classical
results in terms of Bousfield localization. In the third section we discuss
the basic properties of \emph{rigid localizing subcategories} and give a
counterexample of a non-rigid localizing subcategory generated by a set. In
section four we state and prove the classification theorem, the rigid
localizing subcategories in the derived category of quasi-coherent torsion
sheaves on a noetherian formal scheme $\FX$ are in one to one correspondence
with the subsets in the underlying space of $\FX$. The arguments are close in
spirit to \cite{Nct}, with the modifications needed to make them work in the
present context. In the last section we give a description of the acyclization
functor associated to a stable for specialization subset as the derived
functor of the sections with support and connect it to smashing
localizations and to Lipman's idempotent pair. We characterize the localizing
subcategory associated to such a subset by means of homological support.
This result gives us a comparison of Thomason classification and ours for a
noetherian separated scheme. Finally, by adjointness, we obtain also a
description of the localization functor associated to generically stable
subsets. 

The question of describing localizations for subsets that are neither stable
for specialization nor generically stable remains open for the moment.

\begin{ack} We thank Joe Lipman for his patience and interest. His remarks
have allowed us to greatly improve this paper. Also, A.J.L. and L.A.T. wish to
thank Clara Alonso Jerem\'{\i}as for the happy moments she shared with them
during the last part of the preparation of this work
\end {ack}

\section{Basic facts and set-up}

\begin{cosa}\textbf{Preliminaries.} For formal schemes, we will follow the
terminology of \cite[Section 10]{ega1} and of \cite{dfs}. In this paper, we
will always consider noetherian schemes and noetherian formal schemes.

Let $(\FX, \CO_\FX)$ be a noetherian formal scheme and let $\CI$ be an
ideal of definition of $\FX$. In what follows, we will identify an usual
(noetherian) scheme with a formal scheme whose ideal of definition is 0.
Denote by $\CA(\FX)$ the category of all $\CO_\FX$-modules. The powers of
$\CI$ define a torsion class (see \cite[pp. 139-141]{St}) whose associated
torsion functor is 
\[ \varGamma_\CI \CF := \dirlim{n > 0} \shom_{\CO_\FX}(\CO_\FX/\CI^n, \CF)
\] for $\CF \in \CA(\FX)$. This functor does not depend on $\CI$ but on the
topology it determines in the rings of sections of $\CO_\FX$, therefore we
will denote it by $\varGamma'_\FX$. Let $\CA_\ts(\FX)$ be the full
subcategory of $\CA(\FX)$ consisting of sheaves $\CF$ such that
$\varGamma'_\FX \CF = \CF$; it is a \emph{plump} subcategory of $\CA(\FX)$.
This means it is closed for kernels, cokernels and extensions (\emph{cfr.}
\cite[beginning of \S 1]{dfs}). Most important for us is the subcategory 
$\CA_\qct(\FX) : = \CA_\ts(\FX) \cap \CA_\qc(\FX)$. It is again a plump
subcategory of $\CA(\FX)$ by \cite[Corollary 5.1.3]{dfs} and it defines a
triangulated subcategory of $\D(\FX) := \D(\CA(\FX))$, the derived category
of $\CA(\FX)$, it is $\D_\qct(\FX)$, the full subcategory of
$\D(\FX)$ formed by complexes whose homology lies in  $\CA_\qct(\FX)$. If
$\FX = X$ is an usual scheme then $\CA_\ts(X) = \CA(X)$ and 
$\CA_\qct(X) = \CA_\qc(X)$. 

The inclusion functor $\CA_\qct(\FX) \to \CA(\FX)$ has a right adjoint
denoted $Q^\ts_\FX$ (see \cite[Corollary 5.1.5]{dfs}). By the existence
of K-injective resolutions (\cite[Theorem 4.5]{S} or \cite[Theorem
5.4]{AJS}) it is possible to get right-derived functors from functors with
source a category of sheaves, as a consequence we have a functor 
$\R Q^\ts_\FX : \D(\FX) \to \D(\CA_\qct(\FX))$. If $\FX$ is either separated
or of finite Krull dimension, this functor induces an equivalence between
$\D_\qct(\FX)$ and $\D(\CA_\qct(\FX))$ by \cite[Proposition 5.3.1]{dfs}. In
these cases, we will identify $\D(\CA_\qct(\FX))$ and $\D_\qct(\FX)$. To
avoid potential confusions, let us point out that all left and right derived
functors defined over $\D_\qct(\FX)$, or over $\D(\CA_\qct(\FX))$ when this
category is equivalent to the former, are defined using K-flat and K-injective
resolutions in $\K(\FX)$.

The categories $\CA_\qct(\FX)$ and $\CA(\FX)$ are Grothendieck categories so
we can apply the machinery developed in \cite{AJS}. In particular, if $\CL$ is
the smallest localizing subcategory of $\D(\CA_\qct(\FX))$ or of $\D(\FX)$
that contains a given set, then there is a localization functor $\ell$ such
that $\CL$ is the full subcategory of $\D(\CA_\qct(\FX))$ or of $\D(\FX)$,
respectively, whose objects are sent to 0 by $\ell$ (see \cite[Theorem
5.7]{AJS}). The category $\D_\qct(\FX)$ is a localizing subcategory of
$\D(\FX)$, therefore if $\CL$ is the smallest localizing subcategory of
$\D(\FX)$ that contains a given set of objects in $\D_\qct(\FX)$, the
localization functor $\ell$ defined over $\D(\FX)$ lands inside
$\D_\qct(\FX)$, therefore $\CL$ is characterized again as the full subcategory
of $\D_\qct(\FX)$ whose objects are sent to 0 by $\ell$. If $\FX$ is either
separated or of finite Krull dimension, the localizations of $\D_\qct(\FX)$
are identified with those of $\D(\CA_\qct(\FX))$. For the general formalism of
Bousfield localization in triangulated categories the reader may consult
\cite[\S 1]{AJS}.

\end{cosa}
\begin{cosa}\textbf{Monoidal structures.}
The categories $\CA(\FX)$ and $\CA_\qct(\FX)$ are \emph{symmetric closed}, in
the sense of Eilenberg and Kelly, see \cite{lajolla}.  For every $\CF \in
\K(\CA(\FX))$ there is a K-flat resolution $\CP_\CF \to \CF$, this follows
from \cite[Proposition 5.6]{S}. As a consequence, there exists a derived
functor:
\[ \CF \otimes^{\LL}_{\CO_\FX} - : \D(\FX) \to \D(\FX)
\] defined by $\CF \otimes^{\LL}_{\CO_\FX} \CG = \CP_\CF \otimes_{\CO_\FX}
\CG$. Also the functor $\shom^\cdot_{\CO_\FX}(\CF,-)$ has a right derived
functor defined by $\R\shom^\cdot_{\CO_\FX}(\CF,\CG) =
\shom^\cdot_{\CO_\FX}(\CF,\CJ_\CG)$ where $\CG \to \CJ_\CG$ denotes a
K-injective resolution of $\CG$. The usual relations hold providing $\D(\FX)$
with the structure of symmetric closed category. Observe that the unit object
is $\CO_\FX$.

Given $\CF, \CG \in \D_\qct(\FX)$, the complex $\CF \otimes^{\LL}_{\CO_\FX}
\CG$ has quasi-coherent torsion homology. Indeed, it is a local question, and
for affine noetherian formal schemes, a complex in $\D_\qct(\FX)$ is
quasi-isomorphic to a complex made by locally free sheaves so the homology
of $\CF \otimes^{\LL}_{\CO_\FX} \CG$ is quasi-coherent. Furthermore for any
$\CF \in \D_\ts(\FX)$ and $\CE \in \D(\FX)$, the complex
$\CF \otimes^{\LL}_{\CO_\FX} \CE \in \D_\ts(\FX)$. Again, this is a local
question so it can be checked using \cite[Proposition 5.2.1 (a)]{dfs} and the
complex $\CK^\cdot_\infty$ in its proof. Therefore, for each $\CF \in
\D_\qct(\FX)$, the functor
$\CF \otimes^{\LL}_{\CO_\FX}- : \D_\qct(\FX) \to \D(\FX)$ takes values in
$\D_\qct(\FX)$. So it provides an internal tensor product. One can see that
the category $\D_\qct(\FX)$ has a symmetric monoidal structure. The unit
object is $\R\varGamma'_\FX \CO_\FX$ where by $\R\varGamma'_\FX$ we denote
the right-derived functor of $\varGamma'_\FX$. We will denote this object by
$\CO'_\FX$ for convenience.

If furthermore $\FX$ is either separated or of finite Krull dimension, the
category $\D_\qct(\FX) = \D(\CA_\qct(\FX))$ possesses the richer structure of
symmetric closed category. The internal hom is defined as:
\[\dhom^\cdot_{\FX}(\CF,\CG) :=
\R{}Q^\ts_\FX\R\shom^\cdot_{\CO_\FX}(\CF,\CG)
\]  for $\CF, \CG \in \D_\qct(\FX)$. It is also important to note that the
$\otimes$--hom adjunction is internal, \ie it holds replacing the usual
hom-group with the internal hom we have just defined, namely, we have a
canonical isomorphism:
\[ \dhom^\cdot_{\FX}(\CF\otimes^{\LL}_{\CO_\FX}\CG, \CM) =
\dhom^\cdot_{\FX}(\CF,\dhom^\cdot_{\FX}(\CG,\CM))
\]
where $\CF, \CG \text{ and } \CM \in \D_\qct(\FX)$.

If the reader is only interested in usual schemes, then it is enough to
consider the quasi-coherence of the derived tensor product. In this case the
topology in the sections of the structural sheaf is discrete, $\varGamma'_X$
is the identity functor and so the unit object is $\CO_X$. For the internal
hom-sheaf, in the separated or finite Krull dimension case, one uses the
derived ``coherator'' functor $\R{}Q$ defined in \cite[\S 3]{I} taking:
\[\dhom^\cdot_{X}(\CF,\CG) := \R{}Q\R\shom^\cdot_{\CO_X}(\CF,\CG)
\]
for $\CF \text{ and } \CG \in \D(\CA_\qc(X))$.

\end{cosa}

\section{Cohomology with supports on formal schemes}
\begin{cosa}\label{algsup}\textbf{Algebraic supports.}
Given $\CF \in \D_\qct(\FX)$ and $Z \subset \FX$ a closed subset, for the
right derived functor of sheaf of sections with support along $Z$ we have that
$\R{}\varGamma_Z \CF \in \D_\qct(\FX)$ because in the distinguished
triangle
\begin{equation} \label{Zcerrado} 
\R{}\varGamma_Z \CF \to \CF \to \R{}j_*j^*\CF \overset{+}{\to},
\end{equation}
where $j : \FX \setminus Z \inc \FX$ denotes the canonical open embedding,
$\R{}j_*j^*\CF \in \D_\qct(\FX)$
\cite[Proposition 5.2.6 and Corollary 5.2.11]{dfs}. On the other hand, the
closed subset $Z$ is the support of a coherent sheaf $\CO_\FX / \CQ$ where
$\CQ$ is an open coherent ideal in $\CO_\FX$. The functor
\[  \varGamma'_Z := \varGamma_\CQ = 
                   \dirlim{n > 0} \shom_{\CO_\FX}(\CO_\FX/\CQ^n, -)
\] of ``sections with algebraic support along $Z$'' does not depend on $\CQ$
but only on $Z$. The natural map  $\varGamma'_Z \to \varGamma_Z$ is an
isomorphism when applied to sheaves in  $\CA_\qct(\FX)$. Furthermore the
natural morphism in $\D(\FX)$ obtained by deriving 
$\theta_{Z,\CF} : \R\varGamma'_Z \CF \to \R\varGamma_Z \CF$ is an isomorphism
for all $\CF \in \D_\qct(\FX)$. Indeed, this is a local question, so we can
assume that $\FX$ is affine with $\FX = \spf(A)$ where $A$ is a noetherian
adic ring. Let $\kappa : \spf(A) \to \spec(A)$ be the canonical map. Let $X
:= \spec(A)$. The set $Z$ can be considered as a closed subset of either $X$
or $\FX$. We will use $\varGamma'_Z$ and $\varGamma_Z$ for the corresponding
pair of endofunctors in $\CA(\FX)$ and $\CA(X)$. This will not cause any
confusion, because the context will make it clear in which category we are
working. By \cite[Proposition 5.2.4]{dfs} it is enough to show that
$\kappa_*\theta_{Z,\CF}$ is an isomorphism. But this is true because the
diagram
\[\begin{CD}
    \kappa_*\R\varGamma'_Z \CF @>\kappa_*\theta_{Z,\CF}>>
                                            \kappa_*\R\varGamma_Z\CF \\
         @VVV                                          @VVV          \\ 
    \R\varGamma'_Z \kappa_*\CF @>>>      \R\varGamma_Z\kappa_*\CF
\end{CD} \]
commutes and all the unlabeled maps are isomorphisms (for the map in the
bottom use \emph{loc.\!~cit.} and \cite[Corollary 3.2.4]{AJL1}).

Given $\CE, \CF \in \D_{\qct}(\FX)$ there is a bifunctorial map 
\[\psi_{Z}(\CE, \CF) : \CE \otimes^{\LL}_{\CO_\FX} \R\varGamma_Z\CF
\to  \R\varGamma_Z (\CE \otimes^{\LL}_{\CO_\FX} \CF)\]
defined as follows. Assume $\CE$ is K-flat and $\CF$ is K-injective and
choose a quasi-isomorphism $\CE \otimes_{\CO_\FX} \CF \to \CJ$ with
$\CJ$ K-injective. The composed map (of complexes) 
$\CE \otimes_{\CO_\FX} \varGamma_Z\CF \to \CE \otimes_{\CO_\FX} \CF
\to \CJ$ has image into $\varGamma_Z\CJ$ and we define $\psi_{Z}(\CE, \CF)$
to be the resulting factorization
\[\CE \otimes^{\LL}_{\CO_\FX} \R\varGamma_Z\CF 
    \iso \CE \otimes_{\CO_\FX} \varGamma_Z\CF 
    \xto{\psi_{Z}(\CE, \CF)} \varGamma_Z\CJ 
    \iso \R\varGamma_Z(\CE \otimes^{\LL}_{\CO_\FX} \CF).
\]
This map is a quasi-isomorphism if $Z$ is closed. The question is local so
using again \cite[Proposition 5.2.4 and Proposition 5.2.8]{dfs} we restrict 
to the analogous question for an ordinary scheme $X$ and a closed subset $Z
\subset X$. We conclude by \cite[Corollary 3.2.5]{AJL1}.
\end{cosa}

\begin{cosa}\label{ssfs}\textbf{Systems of supports on formal schemes.}
In general, a subset $Z \subset \FX$ \emph{stable for specialization} is a
union $Z = \bigcup_{\alpha \in I}Z_\alpha$ of a directed system of closed
subsets $\{Z_\alpha \,/\, \alpha \in I\}$ of $\FX$ and $\varGamma_Z =
\dirlim{\alpha \in I} \varGamma_{Z_\alpha}$, this corresponds to the
classical case of a ``system of supports''. Writing $\varGamma'_Z =
\dirlim{\alpha \in I} \varGamma'_{Z_\alpha}$ the canonical map
$\varGamma'_Z \to \varGamma_Z$ induces natural maps  $\theta_{Z,\CF}
: \R\varGamma'_Z \CF \to \R\varGamma_Z \CF$  for all $\CF \in \D(\FX)$.
If $\CF \to \CJ$ is a K-injective resolution, we have that
\[ \theta_{Z,\CF} : \R\varGamma'_Z \CF = \varGamma'_Z \CJ = 
\dirlim{\alpha \in I} \varGamma'_{Z_\alpha} \CJ 
\xto{{\smdirlim{\alpha \in I}} \theta_{Z_\alpha,\CF}}
\dirlim{\alpha \in I} \varGamma_{Z_\alpha} \CJ = \R\varGamma_Z \CF
\] therefore, for all $\CF \in \D_\qct(\FX)$, $\theta_{Z,\CF}$ is a
quasi-isomorphism. 

Mimicking the case of a closed subset, for $\CE, \CF \in \D_{\qct}(\FX)$
there is a bifunctorial map 
\[\psi_{Z}(\CE, \CF) : \CE \otimes^{\LL}_{\CO_\FX} \R\varGamma_Z\CF
\to  \R\varGamma_Z (\CE \otimes^{\LL}_{\CO_\FX} \CF)\]
that is a quasi-isomorphism. To check this fact we may assume $\CE$ is K-flat
and $\CF$ is K-injective and choose a quasi-isomorphism 
$\CE \otimes_{\CO_\FX} \CF \to \CJ$ with $\CJ$ a K-injective resolution and
consider the commutativity of the diagram of complexes
\[\begin{CD}
 \CE \otimes^{\LL}_{\CO_\FX} \varGamma_Z \CF
                 @>\psi_{Z}(\CE, \CF)>> \varGamma_Z\CJ \\
     @AAA                         @AAA                 \\ 
 \dirlim{\alpha \in I} (\CE \otimes^{\LL}_{\CO_\FX} \varGamma_{Z_\alpha} \CF)
                 @>{\smdirlim{\alpha \in I}}\psi_{Z_\alpha}(\CE, \CF)>>
                 \dirlim{\alpha \in I}\varGamma_{Z_\alpha} \CJ.
\end{CD} \]
\end{cosa}

\begin{cosa}\label{btfss}\textbf{Bousfield triangles for systems of supports.}
Let $Z \subset \FX$ be a subset stable for specialization as in the previous
paragraph. The endofunctor $\R\varGamma_Z \colon \D_\qct(\FX) \to
\D_\qct(\FX)$ together with the natural transformation
$\rho \colon \R\varGamma_Z \to \id$ is a Bousfield acyclization functor.
Let us see why. We need to check that $\rho$ induces a canonical isomorphism
$\rho (\R\varGamma_Z \CM) = (\R\varGamma_Z \rho)(\CM)$, for all $\CM \in
\D_\qct(\FX)$. Indeed, it follows from the previous paragraph that it is
enough to check this for $\CM \in \D^+_\qct(\FX)$, specifically for $\CM
= \CO'_\FX$. The question is local, so arguing as at the end of 
\ref{algsup}, we can suppose that $\FX = X$ is a noetherian affine scheme and
$\CM$ a bounded-below complex formed by quasi-coherent injective sheaves. In
this case $\varGamma_Z \CM$ is a bounded-below complex formed by
quasi-coherent injective sheaves, too (\emph{cfr.} \cite[Propositions
VI.7.1 and VII.4.5]{St}). But the functor $\varGamma_Z$ is idempotent
from which it follows that 
\[ \R\varGamma_Z\R\varGamma_Z \CM = \varGamma_Z\varGamma_Z \CM =
   \varGamma_Z \CM = \R\varGamma_Z \CM
\]

Using the notation of paragraph \ref{algsup} for a closed subset $Z \subset
\FX$, the triangle (\ref{Zcerrado}) is a Bousfield localization triangle for
each $\CF \in \D_\qct(\FX)$.

In general, let $Z \subset \FX$ be a subset stable for specialization,
therefore it can be considered as the union of a directed system
$\{Z_\alpha \,/\, \alpha \in I\}$ of closed subsets of $\FX$. For every
$\alpha \in I$, let $U_\alpha := \FX \setminus Z_\alpha$ be the complementary
open subset and $j_\alpha\colon U_\alpha \to \FX$ be the canonical open
embedding. Let $L_Z \colon \CA(\FX) \to \CA(\FX)$ be the
endofunctor defined as $L_Z := \dirlim{\alpha \in I} j_{\alpha\,*}
j_\alpha^*$. For every $\CM \in \D_\qct(\FX)$ the triangle
\[ \R\varGamma_Z \CM \overset{\rho(\CM)}{\lto} \CM \lto 
   \R{}L_Z \CM \overset{+}{\lto}
\]
is the Bousfield localization triangle whose associated acyclization functor
is $\R\varGamma_Z$.

For every $\CE \in \D_\qct(\FX)$, the commutative diagram
\[\begin{CD}
    \CE \otimes^{\LL}_{\CO_\FX} \R\varGamma_Z\CO'_\FX 
                     @>\CE \otimes \rho(\CO'_\FX)>>   
                     \CE \otimes^{\LL}_{\CO_\FX} \CO'_\FX   \\
         @V\psi_Z(\CE, \CO'_\FX)V{\wr}V        @VV{\wr}V      \\ 
    \R\varGamma_Z\CE @>\rho(\CE)>>      \CE
\end{CD} \]
can be completed to an isomorphism of distinguished triangles
\[\begin{CD}
    \CE \otimes^{\LL}_{\CO_\FX} \R\varGamma_Z\CO'_\FX
                 @>>>   \CE \otimes^{\LL}_{\CO_\FX} \CO'_\FX  
                 @>>>   \CE \otimes^{\LL}_{\CO_\FX} \R{}L_Z\CO'_\FX
                 @>+>> \\
         @VV{\wr}V            @VV{\wr}V              @VV{\wr}V          \\ 
    \R\varGamma_Z\CE @>>>      \CE @>>> \R{}L_Z\CE @>+>> .
\end{CD} \]

Note that, in particular, $\R\varGamma_Z$ and $\R{}L_Z$ are endofunctors of
$\D_\qct(\FX)$ that commute with coproducts, and two Bousfield
acyclization or localization functors of this kind commute. If $Z, W \subset
\FX$ are stable for specialization subsets then 
$\varGamma_{Z \cap W} = \varGamma_Z \varGamma_W$. One can check, following
the same kind of arguments at the beginning of this subsection, that the
canonical map
$\R\varGamma_{Z \cap W} \CF \to \R\varGamma_Z \R\varGamma_W \CF$ is an
isomorphism for every $\CF \in \D_\qct(\FX)$.
\end{cosa}

\begin{cosa}\label{compL}\textbf{Computing the functor} $\R{}L_{\FX \setminus
\FX_x}$\textbf{.} 
Let $x \in \FX$. Consider the affine formal scheme
$\FX_x := \spf(\widehat{\CO_{\FX,x}})$ where the adic topology in the ring
$\CO_{\FX,x}$ is given by $\CI_x$. If $\FX = \spf{B}$ and $\ip$ is the prime
ideal corresponding to the point $x$, then $\CO_{\FX,x} = B_{\{\ip\}}$.
Denote by $i_x : \FX_x \inc \FX$ the canonical inclusion map. Consider the
functors
\[
\D_\qct(\FX_x)
\underset{i_x^*}{\overset{\R{}i_{x*}}{\rightleftarrows}}
\D_\qct(\FX).
\] which are defined by virtue of \cite[Proposition 5.2.6 and
Corollary 5.2.11]{dfs} using the fact that $i_x$ is an adic map.

Given $\CF_1, \CF_2 \in \D_\qct(\FX)$, we have that
\[
\Hom_{\D(\FX)}(\R\varGamma_{\FX \setminus \FX_x}\CF_1,
\R{}i_{x*}i_x^*\CF_2) \cong
\Hom_{\D(\FX_x)}(i_x^*\R\varGamma_{\FX \setminus \FX_x}\CF_1, i_x^*\CF_2) = 0
\] because $i_x^*\R\varGamma_{\FX \setminus \FX_x}\CF_1 = 0$. Indeed, write 
$\FX \setminus \FX_x = \bigcup_{\alpha \in I}Z_\alpha$ a filtered union of
closed subsets $\{Z_\alpha \,/\, \alpha \in I\}$, and let $\CF_1 \to \CJ$ be a
K-injective resolution then:
\[  i_x^*\R\varGamma_{\FX \setminus \FX_x}\CF_1
  = i_x^*\varGamma_{\FX \setminus \FX_x}\CJ
  = i_x^*\dirlim{\alpha \in I}\varGamma_{Z_\alpha}\CJ
  = \dirlim{\alpha \in I}i_x^*\varGamma_{Z_\alpha}\CJ
  = 0.
\]
It follows that for each $\CF \in \D_\qct(\FX)$ there is a unique map
$\R{}L_{\FX \setminus \FX_x}\CF \to \R{}i_{x*}i_x^*\CF$ making the
following diagram commutative:
\[\begin{CD}
    \CF @>>> \R{}L_{\FX \setminus \FX_x}\CF  \\
    @|       @VV{h_\CF}V                       \\ 
    \CF @>>> \R{}i_{x*}i_x^*\CF.
\end{CD} \]
Furthermore $h$ is a natural transformation of $\Delta$-functors and it is an
isomorphism, \ie $h_\CF$ is a quasi-isomorphism for every $\CF \in
\D_\qct(\FX)$. Let us show this. First of all, we can assume that $\FX$ is
affine. Indeed, choose an affine open subset $\FU \subset \FX$ such that
$x \in \FU$, then one can describe $\FX \setminus \FX_x$  as a filtered union
of closed subsets $\{Z_\alpha \,/\, \alpha \in I\}$ such that each $\FU_\alpha
:= \FX \setminus Z_\alpha$ is an affine open subset of $\FU$. Let us denote
by $j : \FU \inc \FX$, $j_\alpha : \FU_\alpha \inc \FX$ and $i'_x : \FX_x
\inc \FU$ the canonical morphisms. Note that $j \circ i'_x = i_x$. For
every  $\CF \in \D_\qct(\FX)$ we have an isomorphism $\R{}L_{\FX \setminus
\FX_x}\CF \iso \R{}j_*j^*\R{}L_{\FX \setminus \FX_x}\CF$ because
$\R\varGamma_{\FX \setminus \FU}\R{}L_{\FX \setminus\FX_x}\CF =
\R\varGamma_{\FX \setminus \FU}\R\varGamma_{\FX \setminus \FX_x}\R{}L_{\FX
\setminus\FX_x}\CF = 0$ (see \ref{btfss}). Using flat base change
\cite[Proposition 7.2]{dfs} we see that the canonical map $\R{}i_{x*}i_x^*\CF
\iso \R{}j_*j^*\R{}i_{x*}i_x^*\CF$ is also an isomorphism. So, we are left to
prove that $j^*h_\CF$ is an isomorphism, or, equivalently, that $h_{j^*\CF} :
\R{}L_{\FU \setminus \FX_x}(j^*\CF) \to \R{}i'_{x*}{i'}_x^*(j^*\CF)$ is an
isomorphism. Then, let us treat the case $\FX =\spf{A}$ with $A$ a
complete noetherian ring. Both endofunctors $\R{}L_{\FX \setminus \FX_x}$ and
$\R{}i_{x*}i_x^*$ commute with coproducts by \ref{btfss} and
\cite[Proposition 3.5.2]{dfs} respectively. To prove that $h_\CF$ is a
quasi-isomorphism for every 
$\CF \in \D(\CA_\qct(\FX)) = \D_\qct(\FX)$ it is enough to check it for $\CF
\in \CA_\qct(\FX)$, because the smallest localizing subcategory containing
$\CA_\qct(\FX)$ is all of $\D(\CA_\qct(\FX))$. In this case the morphisms
$j_\alpha : U_\alpha \inc \FX$ and  $i_x : \FX_x \inc \FX$ are affine.
Therefore, by  \cite[Lemma 3.4.2]{dfs}, for $\CF \in \CA_\qct(\FX)$ and 
$i > 0$
\begin{align*}
\CH^i(\R{}L_{\FX \setminus \FX_x} \CF) & = 
   \dirlim{\alpha \in I} \CH^i(\R{}j_{\alpha\,*} j^*_\alpha \CF) = 0 \\
\CH^i(\R{}i_{x*}i_x^* \CF) & = 0,
\end{align*}
and for $i = 0$
\[\CH^0(\R{}L_{\FX \setminus \FX_x} \CF) = 
   \dirlim{\alpha}j_{\alpha\,*} j^*_\alpha \CF \overset{\CH^0(h_\CF)}{\lto}
   i_{x*}i_x^* \CF = \CH^0(\R{}i_{x*}i_x^* \CF)
\] is the natural map. Let us show that $\CH^0(h_\CF)$ is an isomorphism.
Using \cite[Proposition 5.2.4]{dfs} we are reduced to the particular case
$\FX = X = \spec{A}$ is an usual affine scheme, $x$ corresponds to a prime
ideal $\ip \subset A$, $M$ is an $A$-module and  $\CF = \widetilde{M}$. Then
$\CH^0(h_\CF)$ corresponds to the canonical isomorphism of $A$-modules 
\[\dirlim{f \in A \setminus \ip}M_f \iso M_\ip.\]

Therefore for $\FX$ a noetherian formal scheme and every $\CF \in
\D_\qct(\FX)$ one has a natural Bousfield triangle
\begin{equation}\label{loctrig}
\R\varGamma_{\FX \setminus \FX_x}\CF \to \CF \to \R{}i_{x*}i_x^*\CF 
  \overset{+}{\to}.
\end{equation}
Recall that the canonical triangle
\[\R\varGamma_{\FX \setminus \FX_x}\CO'_\FX \to \CO'_\FX 
  \to \R{}i_{x*}i_x^*\CO'_\FX \overset{+}{\to}
\] tensored by $\CF$ provides a triangle
\[\R\varGamma_{\FX \setminus \FX_x}\CO'_\FX \otimes^{\LL}_{\CO_\FX} \CF 
  \to \CF \to \R{}i_{x*}i_x^*\CO'_\FX \otimes^{\LL}_{\CO_\FX} \CF
  \overset{+}{\to}
\] that is naturally isomorphic to (\ref{loctrig}) by \ref{btfss}.
\end{cosa}

\section{Rigid localizing subcategories}

Let $\T$ be a triangulated category with all coproducts. This is the case
for $\D(\FX)$ and $\D_\qct(\FX)$ for a noetherian formal scheme $\FX$, and
also for $\D(X)$ and $\D_\qc(X)$ for an usual scheme $X$. A triangulated
subcategory $\CL$ of $\T$ is called \emph{localizing} if it is stable for
coproducts in $\T$. If $\T$ is one of the aforementioned derived categories it
is not ensured that $\CL \subset \T$ is well-behaved with respect to the
tensorial structure. It turns out that we need such compatibility in order to
localize on open subsets. So let us establish the following definition. A
localizing subcategory $\CL \subset \D_\qct(\FX)$ is called \emph{rigid} if
for every $\CF \in \CL$ and $\CG \in \D_\qct(\FX)$, we have that  $\CF
\otimes^{\LL}_{\CO_\FX} \CG \in \CL$. This condition has been independently
considered by Thomason for thick subcategories by the same reason (see
\cite[Definition 3.9]{th}, where they are called
$\otimes$-\emph{subcategories}). Our route to find this condition came from a
paper by one of the authors where localizations are considered in the abelian
context, see \cite[2.3]{saga3}.

\begin{prop} \label{righom} Suppose that $\FX$ is furthermore either separated
or of finite Krull dimension. Let $\CL$ be a localizing subcategory of
$\D(\CA_\qct(\FX))$. If $\CL$ is rigid, then, for every $\CF, \CG \in
\D(\CA_\qct(\FX))$ such that $\CG$ is $\CL$-local (\ie $\CG \in
\CL^{\perp}$), then $\dhom^\cdot_{\FX}(\CF,\CG)$ is $\CL$-local. If moreover
${^\perp{(\CL^\perp)}} = \CL$, the converse is true.
\end{prop}

\begin{proof} Let $\CH \in \CL$, then,
\begin{equation} \label{tenhom} 
   \Hom_{\D(\FX)}(\CH, \dhom^\cdot_{\FX}(\CF,\CG)) =
   \Hom_{\D(\FX)}(\CH \otimes^{\LL}_{\CO_\FX} \CF, \CG) = 0,
\end{equation} 
because $\CG \in \CL^\perp$ and  $\CH \otimes^{\LL}_{\CO_\FX} \CF \in \CL$.
Conversely, if (\ref{tenhom}) holds for every $\CG \in \CL^\perp$, then 
$\CH \otimes^{\LL}_{\CO_\FX} \CF \in {^\perp{(\CL^\perp)}} = \CL$.
\end{proof}

\begin{rem} The condition ${^\perp{(\CL^\perp)}} = \CL$ holds if $\CL$ is
the localizing subcategory of objects whose image is 0 by a Bousfield
localization (see \cite[Proposition 1.6]{AJS}). We will see later
(Corollary \ref{nonsei}) that every rigid localizing subcategory of
$\D(\CA_\qct(\FX))$ arises in this way.
\end{rem}

\begin{prop} \label{affrig} If $\FX$ is affine, every localizing subcategory
of $\D(\CA_\qct(\FX))$ is rigid.
\end{prop}

\begin{proof} Take $\FX = \spf{A}$ where $A$ is a noetherian adic ring.
Every quasi-coherent torsion sheaf comes from an $A$-module and therefore
it has a free resolution. Let $\kappa : \spf{A} \to \spec{A}$ the canonical
morphism and $X := \spec{A}$. Let $\CL$ be a localizing subcategory of
$\D(\CA_\qct(\FX))$. The full subcategory $\T$ of $\D(\CA_\qc(X))$
defined by
\[\T = \{\CN \in \D(\CA_\qc(X)) \,/\, \kappa^*\CN \otimes^{\LL}_{\CO_\FX}
\CM \in \CL, \forall \CM \in \CL \}
\] is triangulated and stable for coproducts. It is clear that $\CO_X \in
\T$, therefore $\T = \D(\CA_\qc(X))$. Now, given $\CG \in
\D(\CA_\qct(\FX))$, $\CG = \kappa^*\kappa_*\CG$ and 
$\kappa_*\CG \in \D(\CA_\qc(X)) = \T$ (\cite[Proposition 5.1.2]{dfs}),
therefore $\CG \otimes^{\LL}_{\CO_\FX} \CM \in \CL$, for every $\CM \in \CL$.
\end{proof}

\begin{ex} Not all localizing subcategories are rigid. Let us show an
example of a non-rigid localizing subcategory. Our example is based in
Thomason's example (\cite[Example 3.13]{th}) of a thick subcategory that it
is not a $\otimes$-subcategory. Consider the projective line over a field
together with its canonical map $\pi \colon \PR^1_k \to \spec{k}$. Denote by
$\D(\PR^1_k)_\pf$ the full subcategory\footnote{Denoted as
$\D(\PR^1_k)_{\text{parf}}$ in \cite{th}.} of $\D(\CA_\qc(\PR^1_k))$ formed by
perfect complexes (\ie quasi-isomorphic to a bounded complex of locally free
finite-type sheaves). Let $\CL$ the smallest localizing subcategory of
$\D(\CA_\qc(\PR^1_k))$ generated by $\CE := \LL\pi^*\widetilde{k}$.
Note that $\CE \in \D(\PR^1_k)_\pf$ and that $\CL$ is the smallest localizing
subcategory that contains the thick subcategory 
$\CA = \{ \CF \in \D(\PR^1_k)_\pf \,/\, \LL\pi^*\R\pi_*\CF = \CF \}$, which
is a thick subcategory of $\D(\PR^1_k)_\pf$, constructed by Thomason in
\emph{loc. cit.} Every object $\CM \in \CL$ is such that $\LL\pi^*\R\pi_*\CM
= \CM$ because both $\LL\pi^*$ and $\R\pi_*$ commute with coproducts and the
equality holds for $\CE$. Observe that $\CL$ is the essential image of
$\D(\CA_\qc(\spec{k}))$ by the functor $\LL\pi^*$. The localizing category
$\CL$ is not rigid. Indeed, take $\CM \in \CL$, $\CM \neq 0$, we will show
that $\CM \otimes \CO(-1) \notin \CL$. Let $\CF := \R\pi_*\CM$, then
\begin{align*}
      \R\pi_*(\CM \otimes \CO(-1))
          & =      \R\pi_*(\LL\pi^*(\CF) \otimes \CO(-1)) 
                                 \tag{\cite[(3.9.4)]{LDC}} \\
          & \simeq \CF \otimes \R\pi_*\CO(-1) \tag{\cite[2.12.16]{ega3}}\\
          & \simeq 0.
\end{align*}
We conclude that $\CM \otimes \CO(-1)$ is not an object in $\CL$ because 
$\CM \otimes \CO(-1) \neq 0 = \LL\pi^*\R\pi_*(\CM \otimes \CO(-1))$.
\end{ex}

\begin{rem} The rigidity condition may seem strange but, in fact, these are
the localizations that behave well when restricted to open subsets and ``are
detected'' by ample sheaves when they exist. We suggest the interested reader
to adapt \cite[Proposition 3.11]{th} and its corollary to our situation.
We will not get into these details because we do not need them.
\end{rem}

\section{Localizing subcategories and subsets}

We keep denoting by $\FX$ a noetherian formal scheme and $\CI$ its
ideal of definition. Let $x \in \FX$, we denote by $i_x : \FX_x \inc \FX$
the canonical inclusion map where $\FX_x = \spf(\widehat{\CO_{\FX,x}})$
(completion with respect to $\CI_x$).
 
We will denote by $\kk(x)$ the residue field of the local ring
$\widehat{\CO_{\FX,x}}$, or, equivalently, of $\CO_{\FX,x}$, by $\CK_x$ the
quasi-coherent torsion sheaf over $\spf(\widehat{\CO_{\FX,x}})$ associated to
the $\widehat{\CO_{\FX,x}}$-module $\kk(x)$ and $\CK(x) := 
\R{}i_{x*}(\CK_x)$. Observe that $\CK(x) =
\R{}\varGamma_{\overline{\{x\}}}\CK(x) = \R{}i_{x*}i_x^*\CK(x)$. If $\FX = X$
is an usual scheme and $x$ is a closed point, $\CK(x)$ has been denoted
$\CO_x$ in recent literature, but we will not use this notation to avoid
potential confusions.

Let $Z$ be any subset of the underlying space of $\FX$. We define the
subcategory $\CL_Z$ as the smallest localizing subcategory of
$\D_\qct(\FX)$ that contains the set of quasi-coherent torsion sheaves
$\{\CK(x) / x \in Z\}$. If $Z = \{x\}$ we will denote $\CL_Z$ simply by
$\CL_x$. Note that if $x \in Z$, then $\CL_x \subset \CL_Z$.

\begin{lem}\label{critpunt} If $\CF \in \D_\qct(\FX)$ and $x \in \FX$,
then $\R{}\varGamma_{\overline{\{x\}}}(\R{}i_{x*}i_x^*\CF)$ belongs to the
localizing subcategory $\CL_x$.
\end{lem}

\begin{proof} Let $\CQ_0$ be a sheaf of coherent ideals in $\CO_\FX$ such
that $\supp(\CO_\FX / \CQ_0) = \overline{\{x\}}$ and denote $\CQ :=
i_x^*\CQ_0$. Recall, by \cite[\S 5.4]{dfs}
\begin{align*}
      \R{}\varGamma_{\overline{\{x\}}}(\R{}i_{x*}i_x^*\CF)
          & = \holim{n > 0} \shom_{\CO_\FX}(\CO_\FX/\CQ_0^n,
                                                   i_{x*}\CJ) \\
          & \cong \holim{n > 0}\R{}i_{x*}\shom_{\CO_{\FX_x}}(
                                       \CO_{\FX_x}/\CQ^n, \CJ) 
\end{align*}
where $i_x^*\CF \to \CJ$ is a K-injective resolution.

Let $\CG := \dirlim{n > 0} \shom_{\CO_{\FX_x}}(\CO_{\FX_x}/\CQ^n, \CJ)$ and
let us consider the filtration
\[0 = \CG_0 \subset \CG_1 \subset \CG_2 \subset \dots \subset \CG
\]
where $\CG_n := \shom_{\CO_{\FX_x}}(\CO_{\FX_x}/\CQ^n, \CJ)$ \ie the
subcomplex of $\CJ$ annihilated by $\CQ^n$. The successive quotients
$\CG_n/\CG_{n-1}$ are complexes of quasi-coherent
$\CK_x$-modules and, therefore, isomorphic in
$\D(\CA_\qct(\FX_x))$ to a direct sum of shifts of $\CK_x$.
The functor $\R{}i_{x*}$ preserves coproducts, therefore every
\(\R{}i_{x*}(\CG_n/\CG_{n-1})\) is an object of $\CL_x$. We deduce by
induction, using the distinguished triangles
\[\R{}i_{x*}\CG_{n-1} \to \R{}i_{x*}\CG_{n} \to 
  \R{}i_{x*}(\CG_n / \CG_{n-1}) \overset{+}{\to}
\] that every $\R{}i_{x*}\CG_n$ is in $\CL_x$ for every $n \in \NN$. But we
have
\[ \R{}\varGamma_{\overline{\{x\}}}(\R{}i_{x*}i_x^*\CF)
    \cong \holim{n > 0} \R{}i_{x*}\CG_n,
\] and the result follows from the fact that a localizing subcategory is
stable for homotopy direct limits \cite[Lemma 3.5 and its proof]{AJS}.
\end{proof} 

Let $E_x$ be an injective hull of the $\CO_{\FX,x}$-module $\kk(x)$, then
$E_x$ is a $\CI_x$-torsion $\widehat{\CO_{\FX,x}}$-module. Let then $\CE_x$ be
the sheaf in $\CA_\qct(\FX_x)$ determined by $\Gamma(\FX_x,\CE_x) = E_x$

\begin{cor}\label{capiny} The object 
$\CE(x) := \R{}i_{x*}\CE_x$ belongs to $\CL_x$. 
\end{cor}

\begin{proof} Use the previous lemma and the fact that
$\CE(x) = \R{}\varGamma_{\overline{\{x\}}}(\R{}i_{x*}i_x^*\CE(x))$.
\end{proof} 

\begin{lem}\label{loctens} Let $\CM \in \D_\qct(\FX)$ and $\CL$ the
smallest localizing subcategory of $\D_\qct(\FX)$ that contains $\CM$.
If $\CG \in \D_\qct(\FX)$ is such that $\CM \otimes^{\LL}_{\CO_\FX} \CG = 0$
then $\CF \otimes^{\LL}_{\CO_\FX} \CG = 0$, for every $\CF \in \CL$.
\end{lem}

\begin{proof} The $\Delta$-functor $- \otimes^{\LL}_{\CO_\FX} \CG$
preserves coproducts and therefore the full subcategory whose objects are
those $\CF \in \CL$ such that $\CF \otimes^{\LL}_{\CO_\FX} \CG = 0$ is
localizing, but it contains $\CM$, therefore it is $\CL$.
\end{proof} 

\begin{prop}\label{smalloc} The smallest localizing subcategory $\CL$ of
$\D_\qct(\FX)$ that contains $\CK(x)$ for every $x \in \FX$ is the
whole $\D_\qct(\FX)$.
\end{prop}

\begin{proof} Let $\CF \in \D_\qct(\FX)$ and  $\CC$ denote the family of
subsets $Y \subset \FX$ stable for specialization such that 
$\R{}\varGamma_Y \CF \in \CL$. If $\{W_\alpha\}_{\alpha \in I}$ is a chain in
$\CC$ then
\[\R{}\varGamma_{\cup W_\alpha} \CF = 
\dirlim{\alpha \in I} \varGamma_{W_\alpha} \CJ,\]
for a K-injective resolution $\CF \to \CJ$. By \cite[Theorem 2.2 and
Theorem 3.1]{AJS} 
$\R{}\varGamma_{\cup W_\alpha} \CF = \varGamma_{\cup W_\alpha} \CJ \in \CL$,
because each $\R{}\varGamma_{W_\alpha} \CF = \varGamma_{W_\alpha} \CJ \in
\CL$, so $\cup W_\alpha \in \CC$.

The set $\CC$ is stable for filtered unions, therefore, there is a maximal
element in $\CC$ which we will denote by $W$. We will see that $W = \FX$ from
which it follows that $\CF \cong \R{}\varGamma_\FX \CF \in \CL$. 

Indeed, otherwise suppose $\FX \setminus W \neq \emptyset$. As $\FX$ is
noetherian the family of closed subsets
\[ \CC' = \{\overline{\{z\}} / z\in \FX \text{ and } 
\overline{\{z\}} \cap (\FX \setminus W) \neq \emptyset \}
\] has a minimal subset $\overline{\{y\}}$. If $x \in \overline{\{y\}} \cap
(\FX \setminus W)$, then $\overline{\{x\}} \in \CC'$, but
$\overline{\{y\}}$ is minimal, so $x = y$ and $W \cup \overline{\{y\}} = W
\cup \{y\}$. Consider now the inclusion
$i_y : \FX_y \to \FX$ and the distinguished triangle in $\D_\qct(\FX)$
\[ \R{}\varGamma_W \CF \lto \R{}\varGamma_{W \cup \{y\}} \CF 
\lto \R{}\varGamma_{\overline{\{y\}}}(\R{}i_{y\,*}i_y^*\CF) \overset{+}{\lto}
\] obtained applying $\R{}\varGamma_{W \cup \{y\}}$ to the
canonical triangle
\[\R{}\varGamma_{\FX \setminus \FX_y} \CF \lto \CF 
\lto \R{}i_{y\,*}i_y^*\CF \overset{+}{\lto}.
\]
We deduce that $\R{}\varGamma_{W \cup \{y\}} \CF \in \CL$, because 
$W \in \CC$ and $\R{}\varGamma_{\overline{\{y\}}}(\R{}i_{y\,*}i_y^*\CF) \in
\CL_y \subset \CL$ by Lemma \ref{critpunt}, contradicting the maximality of
$W$.
\end{proof} 

\begin{cor} Let $\CG \in \D_\qct(\FX)$. We have that $\CG = 0$ if, and
only if, $\Hom_{\D(\FX)}(\CK(x)[n],\CG) = 0$ for all $x \in \FX$ and $n \in
\ZZ$.
\end{cor}

\begin{proof} Immediate from Proposition \ref{smalloc}.
\end{proof}

\begin{cor}\label{crithazK} Let $\CG \in \D_\qct(\FX)$ be such that
$\CK(x) \otimes^{\LL}_{\CO_\FX} \CG = 0$ for every $x \in \FX$, then $\CG =0$.
\end{cor}

\begin{proof} It is a consequence of Proposition \ref{smalloc} and Lemma
\ref{loctens}.
\end{proof}

\begin{lem} \label{Kortog} 
If $x \neq y$, then $\CK(x) \otimes^{\LL}_{\CO_\FX} \CK(y) = 0$.
\end{lem}

\begin{proof} There exist an affine open subset $\FU \subset \FX$ such that
it only contains one of the points, for instance assume that
 $x \in \FU$ and $y \notin \FU$. Denote by $j : \FU
\inc \FX$ the canonical inclusion map. Now, using \ref{compL}
\begin{align*}
      \CK(x) \otimes^{\LL}_{\CO_\FX} \CK(y)
          & \cong \R{}j_*j^*\CK(x) \otimes^{\LL}_{\CO_\FX} \CK(y)    \\
          & \cong \R{}j_*j^*\CO'_\FX \otimes^{\LL}_{\CO_\FX} \CK(x)
                                      \otimes^{\LL}_{\CO_\FX} \CK(y) \\
          & \cong \CK(x) \otimes^{\LL}_{\CO_\FX} \R{}j_*j^*\CK(y)    \\
          & =      0
\end{align*}
because $j^*\CK(y) = 0$.
\end{proof}

\begin{cor}\label{critrig} For every subset $Z \subset \FX$, the localizing
subcategory $\CL_Z$ is rigid.
\end{cor}

\begin{proof} The full subcategory \(\CS \subset \D_\qct(\FX)\) defined
by 
\[\CS = \{\CN \in \D_\qct(\FX) \,/\, 
\CN \otimes^{\LL}_{\CO_\FX} \CM \in \CL_Z,\, \forall \CM \in \CL_Z\}\]
is a localizing subcategory of $\D_\qct(\FX)$. For $x \in \FX$, 
$\CK(x) \cong \R{}\varGamma_{\overline{\{x\}}}\R{}i_{x*}i_x^*\CK(x)$,
so using \ref{btfss} and \ref{compL} we have that  
\[\CK(x) \otimes^{\LL}_{\CO_\FX} \CM \cong
\R{}\varGamma_{\overline{\{x\}}}\R{}i_{x*}i_x^*\CK(x) \otimes^{\LL}_{\CO_\FX}
\CM \cong \R{}\varGamma_{\overline{\{x\}}}\R{}i_{x*}i_x^*(\CK(x)
\otimes^{\LL}_{\CO_\FX} \CM).\]
Therefore if $x \in Z$ then 
$\CK(x) \otimes^{\LL}_{\CO_\FX} \CM \in \CL_x \subset \CL_Z$ by Lemma
\ref{critpunt}, and for $x \notin Z$, by Lemma
\ref{Kortog} and Lemma \ref{loctens}, $\CK(x) \otimes^{\LL}_{\CO_\FX} \CM =
0$ it is also in $\CL_Z$. Necessarily $\CS = \D_\qct(\FX)$ by
Proposition \ref{smalloc}.
\end{proof}
 
\begin{cor} \label{conganma} If $Z$ and $Y$ are subsets of $\FX$ such that
$Z \cap Y = \emptyset$, then $\CF \otimes^{\LL}_{\CO_\FX} \CG = 0$ for every
$\CF \in \CL_Z$ and $\CG \in \CL_Y$.
\end{cor}

\begin{proof} This follows from the previous lemma and Lemma \ref{loctens}.
\end{proof}

\begin{cor} Given $x \in \FX$ and $\CF \in \CL_x$ we have
that\[ \CF = 0 \dimp \CF \otimes^{\LL}_{\CO_\FX} \CK(x) = 0.\]
\end{cor}

\begin{proof} By Lemma \ref{Kortog} and Lemma \ref{loctens}, given $\CF \in
\CL_x$, for all  $y \in \FX$, with $y \neq x$ we have that 
$\CF \otimes^{\LL}_{\CO_\FX} \CK(y) = 0$, therefore if also
$\CF \otimes^{\LL}_{\CO_\FX} \CK(x) = 0$, it follows that
$\CF =  0$ by Corollary \ref{crithazK}.
\end{proof}

\begin{cor} \label{crithaz}  Let $\CL$ be a localizing subcategory of
$\D_\qct(\FX)$ and $\CF \in \D_\qct(\FX)$. If 
$\CK(x) \otimes^{\LL}_{\CO_\FX} \CF \in \CL$ for every $x \in \FX$, then 
$\CF \in \CL$. 
\end{cor}

\begin{proof} Let $\CL' = \{\CG \in \D_\qct(\FX) \,/\, \CG
\otimes^{\LL}_{\CO_\FX} \CF \in \CL\}$. The subcategory $\CL'$ is a
localizing subcategory of $\D_\qct(\FX)$ such that $\CK(x) \in \CL'$ for
all $x \in \FX$. By Proposition \ref{smalloc}, we deduce that $\CL' =
\D_\qct(\FX)$, in particular  
$\CO'_\FX\otimes^{\LL}_{\CO_\FX} \CF = \CF \in \CL$.
\end{proof}

\begin{rem} If the localizing subcategory $\CL$ is rigid then: $\CK(x)
\otimes^{\LL}_{\CO_\FX} \CF \in \CL$ for all $x \in \FX$ if, and only if,
$\CF \in \CL$.
\end{rem}

\begin{thm}\label{ThclasB} For a noetherian formal scheme $\FX$ there is a
bijection between the class of rigid localizing subcategories of
$\D_\qct(\FX)$ and the set of all subsets of $\FX$.
\end{thm}

\begin{proof} Denote by $\loc{(\D_\qct(\FX))}$ the class of rigid
localizing subcategories of
$\D_\qct(\FX)$ and by $\Pa(\FX)$ the set of all subsets of $\FX$. Let
us define a couple of maps:
\[
\loc{(\D_\qct(\FX))}
\underset{\phi}{\overset{\psi}{\rightleftarrows}}
\Pa(\FX)
\] and check that they are mutual inverses.
Define for $Z \subset \FX$, $\phi(Z) := \CL_Z$ which is rigid by Corollary
\ref{critrig}, and for a rigid localizing subcategory $\CL$ of
$\D_\qct(\FX)$, 
$\psi(\CL) := \{x \in \FX  / \exists \CG \in \CL \text{ with } \CK(x)
\otimes^{\LL}_{\CO_\FX} \CG \neq 0\}$.

Let us check first that $\psi \circ \phi = \id$. Let $Z \subset \FX$ and 
$x \in Z$, by definition $\CK(x) \in \CL_Z$ and clearly
$\CK(x) \otimes^{\LL}_{\CO_\FX} \CK(x) \neq 0$ by Corollary \ref{crithazK} 
and Lemma \ref{Kortog}, therefore $x \in \psi(\phi(Z))$, so $Z \subset
\psi(\phi(Z))$. Conversely let $x \in \psi(\phi(Z))$, by definition there is
$\CG \in \CL_Z$ such that $\CK(x) \otimes^{\LL}_{\CO_\FX} \CG \neq 0$, by
Corollary \ref{conganma}, $x \in Z$.

Now we have to prove that $\phi \circ \psi = \id$. Let $\CL$ be a rigid
localizing subcategory of $\D_\qct(\FX)$. We will see first that
$\CL_{\psi(\CL)} \subset \CL$ and for this it will be enough to check that
$\CK(x) \in \CL$ for every $x \in \psi(\CL)$. So let $x \in \psi(\CL)$,
there is a $\CG \in \CL$ such that 
$\CK(x) \otimes^{\LL}_{\CO_\FX} \CG \neq 0$. On the other hand 
$\CK(x) \otimes^{\LL}_{\CO_\FX} \CG$ belongs to $\CL$ because $\CL$ is rigid.
We have that
\[ \CK(x) \otimes^{\LL}_{\CO_\FX} \CG \cong \bigoplus_{\alpha \in
J}\CF_\alpha
\] where $J$ is a set of indices and $\CF_\alpha = \CK(x)[s_\alpha]$ with
$s_\alpha \in \ZZ$. Indeed, it is enough to take a free resolution 
$\CM \to i^*_x \CG$ of the complex of quasi-coherent torsion
$\CO_{\FX_x}$-modules  $i^*_x \CG$ and to consider the chain of natural
isomorphisms
\begin{align*}
    \CK(x) \otimes^{\LL}_{\CO_\FX} \CG
        & \cong \R{}i_{x*} i^*_x(\CK(x) \otimes^{\LL}_{\CO_\FX} \CG) \\
        & \cong \R{}i_{x*}(\CK_x \otimes^{\LL}_{\CO_{\FX_x}} i^*_x \CG)
             \tag{\cite[(3.2.4)]{LDC}}\\
        & \cong \R{}i_{x*}(\CK_x \otimes^{\LL}_{\CO_{\FX_x}} \CM)
\end{align*}
and use the fact that both functors $\CK_x\otimes^{\LL}_{\CO_{\FX_x}}-$
and $\R{}i_{x*}$ commute with coproducts. But $\CL$ is localizing, so stable
for coproducts and, as a consequence, for direct summands (see \cite{BN}\/
or \cite[footnote, p.~227]{AJS}). From this, 
$\bigoplus_{\alpha \in J}\CF_\alpha \in \CL$ implies $\CK(x) \in \CL$, as
required. Finally, let us see that $\CL \subset \CL_{\psi(\CL)}$. Let $\CF \in
\CL$, by Corollary \ref{crithaz} to see that $\CF \in \CL_{\psi(\CL)}$ it is
enough to prove that $\CK(x) \otimes^{\LL}_{\CO_\FX} \CF \in \CL_{\psi(\CL)}$
for every $x \in \FX$.  Suppose that the non-trivial situation $\CK(x)
\otimes^{\LL}_{\CO_\FX} \CF \neq 0$ holds. In this case, $x \in \psi(\CL)$,
therefore we conclude that
$\CK(x)
\otimes^{\LL}_{\CO_\FX} \CF \in \CL_x \subset \CL_{\psi(\CL)}$ using Corollary
\ref{critrig} that tells us that $\CK(x) \otimes^{\LL}_{\CO_\FX} \CF$ belongs
to the localizing subcategory generated by $\CK(x)$.
\end{proof}

\begin{rem} In view of Proposition \ref{affrig}, the previous result is a
generalization of \cite[Theorem 2.8]{Nct} from noetherian \emph{affine}
schemes to the bigger category of noetherian formal schemes.
\end{rem}

\begin{cor} 
For a noetherian scheme $X$ there is a bijection between the class
of rigid localizing subcategories of $\D_\qc(X)$ and the set of all subsets
of $X$.
\end{cor}


\begin{cor} \label{nonsei} Every rigid localizing subcategory of
$\D_\qct(\FX)$ has associated a localization functor.
\end{cor}

\begin{proof} Theorem \ref{ThclasB} says that a rigid localizing
subcategory $\CL \subset \D_\qct(\FX)$ is the smallest localizing
subcategory that contains the set $\{\CK(x) \,/\, x \in \psi(\CL)\}$. It
follows from \cite[Theorem 5.7]{AJS} that there is an associated
localization functor for $\CL$.
\end{proof} 

The following consequences of the previous discussion will be used in the
next section.

\begin{lem}\label{locKhom} Let $\CL$ be a rigid localizing subcategory of
$\D_\qct(\FX)$ and $z \in \FX$. If $z \notin \psi(\CL)$, then $\CK(z)$
is a $\CL$-local object.
\end{lem}

\begin{proof} 
Let $\CN \in \D_\qct(\FX)$ consider the natural map
\[ \Hom_{\D(\FX)}(\CN, \CK(z)) \overset{\alpha}{\lto}
   \Hom_{\D(\FX)}(\CN \otimes^{\LL}_{\CO_\FX} \CK(z),
                  \CK(z) \otimes^{\LL}_{\CO_\FX} \CK(z)),
\]
and the map
\[   \Hom_{\D(\FX)}(\CN \otimes^{\LL}_{\CO_\FX} \CK(z),
               \CK(z) \otimes^{\LL}_{\CO_\FX} \CK(z)) 
\overset{\beta}{\lto} \Hom_{\D(\FX)}(\CN, \CK(z))\]
induced by the canonical maps 
\[ \CO_\FX \to \CK(z) \text{\quad and \quad} 
   \CK(z) \otimes^{\LL}_{\CO_\FX} \CK(z) \to \CK(z).
\]
It is clear that $\beta \circ \alpha = \id$. By Corollary \ref{conganma} we
have that $\CN \otimes^{\LL}_{\CO_\FX} \CG = 0$ for all $\CN \in \CL$ and
$\CG \in \CL_z$, and necessarily,
\[ \Hom_{\FX}(\CN, \CK(z)) = 0,
\] therefore, $\CK(z)$ is $\CL$-local.
\end{proof}

\begin{lem}\label{locKhom2} Suppose that $\FX$ is either separated or
of finite Krull dimension and let $\CL$ be a rigid localizing subcategory of
$\D(\CA_\qct(\FX))$ and $z \in \FX$. If $z \notin \psi(\CL)$, then
$\dhom^\cdot_{\FX}(\CG,\CF)$ is a  $\CL$-local objects for every  $\CF \in
\D(\CA_\qct(\FX))$ and $\CG \in \CL_z$.
\end{lem}

\begin{proof}
By Corollary \ref{conganma} we have that 
\[ \Hom_{\D(\FX)}(\CN, \dhom^\cdot_{\FX}(\CG,\CF)) \cong
   \Hom_{\D(\FX)}(\CN \otimes^{\LL}_{\CO_\FX} \CG, \CF) = 0
\] for every $\CN \in \CL$, from which it follows that
$\dhom^\cdot_{\FX}(\CG,\CF)$ is $\CL$-local.
\end{proof}

\section{Compatibility of localization with the monoidal structure}

In this section $\FX$ will denote a noetherian scheme that is either separated
or of finite Krull dimension. Let $\CL$ be a localizing subcategory of
$\D(\CA_\qct(\FX))$ with associated Bousfield localization functor $\ell$.
For every $\CF \in \D(\CA_\qct(\FX))$ there is a canonical distinguished
triangle:
\begin{equation}\label{boustrig}
\gmm \CF \lto \CF \lto \ell \CF \overset{+}{\lto}
\end{equation}
such that $\gmm \CF \in \CL$ and $\ell \CF \in \CL^{\perp}$ (in
other words, $\ell \CF$ is $\CL$-local). The functor $\gmm$ is called the
acyclization or colocalization associated to $\CL$ and was denoted $\ell^a$ in
\cite{AJS}. Here we have changed the notation for clarity. The endofunctors
$\gmm$ and $\ell$ are idempotent in a functorial sense as explained in \S 1
of \emph{loc.~cit.} For all $\CF, \CG \in \D(\CA_\qct(\FX))$ we have the
following canonical isomorphisms
\[ \Hom_{\D(\FX)}(\gmm \CF, \gmm \CG) \liso \Hom_{\D(\FX)}(\gmm \CF, \CG)
\]
\[ \Hom_{\D(\FX)}(\ell \CF, \ell \CG) \liso \Hom_{\D(\FX)}(\CF, \ell \CG)
\] induced by $\gmm \CG \to \CG$ and $\CF \to \ell \CF$, respectively.

\begin{lem} \label{rigcar} With the previous notation, the following are
equivalent
    \begin{enumerate}
        \item The localizing subcategory $\CL$ is rigid.
        \item The natural transformation $\gmm \CG \to \CG$ induces
              isomorphisms
              \[\dhom^\cdot_{\FX}(\gmm \CF, \gmm \CG) \cong
                \dhom^\cdot_{\FX}(\gmm \CF, \CG)
              \] for every $\CF, \CG \in \D(\CA_\qct(\FX))$.
        \item The natural transformation $\CF \to \ell \CF$ induces
              isomorphisms
              \[\dhom^\cdot_{\FX}(\ell \CF, \ell \CG) \cong
                \dhom^\cdot_{\FX}(\CF, \ell \CG)
              \] for every $\CF, \CG \in \D(\CA_\qct(\FX))$.
    \end{enumerate}
\end{lem}

\begin{proof} Let us show (\textit{i}) $\imp$ (\textit{ii}). Let $\CN \in
\D(\CA_\qct(\FX))$, we have the following chain of isomorphisms
\begin{align*}
      \Hom_{\D(\FX)}(\CN, \dhom^\cdot_{\FX}(\gmm \CF, \gmm \CG))
          & \cong 
           \Hom_{\D(\FX)}(\CN \otimes^{\LL}_{\CO_\FX} \gmm \CF, \gmm \CG)\\
          & \overset{a}{\cong}
           \Hom_{\D(\FX)}(\CN \otimes^{\LL}_{\CO_\FX} \gmm \CF, \CG) \\
          & \cong  
           \Hom_{\D(\FX)}(\CN, \dhom^\cdot_{\FX}(\gmm \CF, \CG)) ;
\end{align*}
where $a$ is an isomorphism because $\CL$ is rigid and therefore
$\CN \otimes^{\LL}_{\CO_\FX} \gmm \CF = 
\gmm(\CN \otimes^{\LL}_{\CO_\FX} \gmm \CF)$. Having an isomorphism for every
$\CN \in \D(\CA_\qct(\FX))$ forces the target complexes to be isomorphic. 

We will see now (\textit{ii}) $\imp$ (\textit{iii}). From (\ref{boustrig}), we
have a distinguished triangle
\begin{small}
\[ \dhom^\cdot_{\FX}(\ell \CF, \ell \CG)  \lto 
   \dhom^\cdot_{\FX}(\CF, \ell \CG)       \lto
   \dhom^\cdot_{\FX}(\gmm \CF, \ell \CG)  \overset{+}{\lto}
\] \end{small}
but its third point is null, considering 
\[ \dhom^\cdot_{\FX}(\gmm \CF, \ell \CG) \overset{(ii)}{\cong}
                \dhom^\cdot_{\FX}(\gmm \CF, \gmm \ell \CG) = 0,\]
because $\gmm \ell \CG = 0$.

Finally, let us see that (\textit{iii}) $\imp$ (\textit{i}). Take $\CF \in
\CL$ and  $\CN \in \D(\CA_\qct(\FX))$. To see that 
$\CN \otimes^{\LL}_{\CO_\FX} \CF \in \CL$ it is enough to check that 
$\Hom_{\D(\FX)}(\CN \otimes^{\LL}_{\CO_\FX} \CF, \CG) = 0$ for every $\CG
\in \CL^\perp$ because
${^\perp{(\CL^\perp)}} = \CL$. But this is true:
\begin{equation*}
   \begin{split}
      \Hom_{\D(\FX)}(\CN \otimes^{\LL}_{\CO_\FX} \CF, \CG)
          & \cong 
              \Hom_{\D(\FX)}(\CN, \dhom^\cdot_{\FX}(\CF, \CG)) \\
          & \overset{b}{\cong}
              \Hom_{\D(\FX)}(\CN, \dhom^\cdot_{\FX}(\ell\CF, \ell\CG)) \\
          & = 0,
   \end{split}
\end{equation*}
where $b$ is an isomorphism, as follows from (\textit{iii}) and the fact
that $\CG = \ell\CG$, and the last equality holds because $\CF \in \CL$ and
so $\ell\CF = 0$.
\end{proof}

\begin{ex} \label{sfs} Let $Z$ be a closed subset of $\FX$, or more
generally, a  set stable for specialization\footnote{See \ref{ssfs}.}.
Recall the functor sections with support $\varGamma_Z :
\CA_\qct(\FX) \to \CA_\qct(\FX)$. From paragraph \ref{btfss}, we see that 
$\R{}\varGamma_Z : \D(\CA_\qct(\FX)) \to \D(\CA_\qct(\FX))$, its
derived functor, together with the natural transformation 
$\R{}\varGamma_Z \to \id$ posses the formal properties of an
acyclization such that the associated localizing subcategory
\[\CL  = \{ \CM \in \D(\CA_\qct(\FX)) \,/\, \R{}\varGamma_Z(\CM) = \CM \}\]
is rigid. 

The functor $\R{}\varGamma_Z$ has the following property:
\begin{equation*}
\R{}\varGamma_Z(\CK(x)) =
    \begin{cases}
        0                 & \text{ if } x \notin Z \\
        \CK(x)            & \text{ if } x \in Z .  \\ 
    \end{cases}
\end{equation*}
Indeed, if $x \notin Z$ by Lemma \ref{locKhom}, $\R{}\varGamma_Z(\CK(x))
=0$. On the contrary, if $x \in Z$ then $\CK(x) \in \CL_x \subset \CL_Z$, so
$\R{}\varGamma_Z(\CK(x)) = \CK(x)$. It follows that $\CL$ has to agree
with $\CL_Z$ by Theorem \ref{ThclasB} and, consequently, $\R{}\varGamma_Z$ is
$\gmm_{\!Z}$, the acyclization functor associated to the localizing
subcategory $\CL_Z$. This acyclization functor satisfies a special property,
namely, $\gmm_{\!Z}(\CF \otimes^{\LL}_{\CO_\FX}
\CG)$ and $\CF \otimes^{\LL}_{\CO_\FX} \gmm_{\!Z} \CG$ are canonically
isomorphic, see paragraph \ref{ssfs}.
\end{ex} 

\begin{cosa} Let $\CL$ be a rigid localizing subcategory of
$\D(\CA_\qct(\FX))$ and $\CF, \CG \in \D(\CA_\qct(\FX))$. The morphism
$\CF \otimes^{\LL}_{\CO_\FX} \gmm \CG \to 
\CF \otimes^{\LL}_{\CO_\FX} \CG$ induced by $\gmm \CG \to  \CG$ factors
naturally through $\gmm(\CF \otimes^{\LL}_{\CO_\FX} \CG)$ giving a
natural morphism
\[ t: \CF \otimes^{\LL}_{\CO_\FX} \gmm \CG \lto
\gmm(\CF \otimes^{\LL}_{\CO_\FX} \CG). 
\] Let us denote by
\[ p: \CF \otimes^{\LL}_{\CO_\FX} \ell \CG \lto
\ell(\CF \otimes^{\LL}_{\CO_\FX} \CG) 
\] a morphism such that the diagram
\[\begin{CD}
\CF\otimes^{\LL}_{\CO_\FX}\gmm\CG   @>>> \CF\otimes^{\LL}_{\CO_\FX}\CG
             @>>> \CF\otimes^{\LL}_{\CO_\FX}\ell\CG @>+>>     \\ 
             @VVtV                @|             @VVpV             \\ 
\gmm(\CF\otimes^{\LL}_{\CO_\FX}\CG) @>>> \CF\otimes^{\LL}_{\CO_\FX}\CG
             @>>> \ell(\CF\otimes^{\LL}_{\CO_\FX}\CG) @>+>>
\end{CD}
\] is a morphism of distinguished triangles. In fact, the triangle is
functorial in the sense that the map $p$ is uniquely determined by $t$
due to the fact that $\Hom_{\D(\FX)}(\CF\otimes^{\LL}_{\CO_\FX}\gmm\CG,
\ell(\CF\otimes^{\LL}_{\CO_\FX}\CG)[-1]) = 0$.

We say that the localization
$\ell$ is $\otimes$\emph{-compatible} (or that $\CL$ is $\otimes$-compatible
or that $\gmm$ is $\otimes$-compatible) if the canonical morphism $t$, or
equivalently $p$, is an isomorphism.
\end{cosa}

We remind the reader our convention that $\CO'_\FX$ denotes
$\R\varGamma'_\FX \CO_\FX$.

\begin{thm}\label{tenscomp}In the previous hypothesis we have the following
equivalent statements:
    \begin{enumerate}
        \item \label{sm1} The localization associated to $\CL$ is  
                          $\otimes$-compatible.
        \item \label{sm1'} For every $\CE  \in \CL^{\perp}$ and
                          $\CF \in \D(\CA_\qct(\FX))$ we have that   
                          $\CF\otimes^{\LL}_{\CO_\FX}\CE \in \CL^{\perp}$.
        \item \label{sm2} The functor $\ell$ preserves coproducts.
        \item \label{sm3} A coproduct of $\CL$-local objects is $\CL$-local.
        \item \label{sm4} The set $Z := \psi(\CL)$ is stable for
                          specialization and its associated acyclization
                          functor is $\gmm = \R{}\varGamma_Z$.
    \end{enumerate}
\end{thm}

\begin{proof}
Let us begin proving the non-trivial part of (\textit{\ref{sm1}}) $\dimp$
(\textit{\ref{sm1'}}). Indeed, suppose that (\textit{\ref{sm1'}}) holds and
for $\CF$, $\CG \in \D(\CA_\qct(\FX))$ consider the triangle
\[ \CF\otimes^{\LL}_{\CO_\FX}\gmm\CG \lto \CF\otimes^{\LL}_{\CO_\FX}\CG
           \lto \CF\otimes^{\LL}_{\CO_\FX}\ell\CG \overset{+}{\lto}
\] we have that $\CF\otimes^{\LL}_{\CO_\FX}\gmm\CG \in \CL$ because $\CL$ is
rigid, on the other hand $\CF\otimes^{\LL}_{\CO_\FX}\ell\CG \in \CL^{\perp}$
because $\ell\CG \in \CL^{\perp}$. The fact that the natural maps
\[\CF\otimes^{\LL}_{\CO_\FX}\gmm\CG \overset{t}{\lto}
  \gmm(\CF\otimes^{\LL}_{\CO_\FX}\CG)
\text{ and }
  \CF\otimes^{\LL}_{\CO_\FX}\ell\CG \overset{p}{\lto}
  \ell(\CF\otimes^{\LL}_{\CO_\FX}\CG)\] are isomorphisms follow from
\cite[Proposition 1.6, (\textit{vi}) $\imp$ (\textit{i})]{AJS}.

Let us see now that (\textit{\ref{sm1}}) $\imp$
(\textit{\ref{sm2}}). If the localization associated to $\CL$ is
$\otimes$-compatible we have that, for $\CF \in \D(\CA_\qct(\FX))$,
\[ \ell \CF \cong \CF \otimes^{\LL}_{\CO_\FX} \ell \CO'_\FX
\] from which is clear that $\ell$ preserves coproducts.

The implication (\textit{\ref{sm2}}) $\imp$ (\textit{\ref{sm3}}) is
obvious because $\ell \CF \cong \CF$ if, and only if, $\CF \in \CL^{\perp}$.
To see that (\textit{\ref{sm3}}) $\imp$ (\textit{\ref{sm4}}), we will use an
argument similar to the one in the affine case (\cite[Lemma 3.7]{Nct}).
Assume that $\ell$ preserves coproducts. Let $x \in Z$ and $z \in
\overline{\{x\}}$. If $z \notin Z$, then $\CK(z) \in \CL^{\perp}$ (Lemma
\ref{locKhom}) and
$\CL_z \subset \CL^{\perp}$ because by (\textit{\ref{sm3}}) $\CL^{\perp}$
is localizing, and it follows by Corollary \ref{capiny} that also
$\CE(z) \in \CL^{\perp}$. But $\CE(x) \in \CL$ which contradicts the
existence on a non-zero map $\CE(x) \to \CE(z)$ because $z \in
\overline{\{x\}}$. Therefore $Z$ is stable for specialization and $\gmm
\cong \R{}\varGamma_Z$ by the example on page \pageref{sfs}. The same example
shows that (\textit{\ref{sm4}}) $\imp$ (\textit{\ref{sm1}}).
\end{proof}

\begin{rem} In the category of stable homotopy, $\st$, the localizations for
which condition (\textit{\ref{sm2}}) holds are called \emph{smashing}.
This can be characterized by a condition analogous to (\textit{\ref{sm1}})
in terms of its monoidal structure via the smash product, $\wedge$. So, the
previous result classifies smashing localizations in $\D(\CA_\qct(\FX))$.
\end{rem} 

\begin{cor} \label{clastens} 
There is a bijection between the class of $\otimes$-compatible localizations
of $\D(\CA_\qct(\FX))$ and the set of subsets stable for specialization of
$\FX$.
\end{cor}

In \cite[\S 1.4]{Lrtfs}, Lipman defines an \emph{idempotent pair} for a closed
category. In the case in which the closed category is $\D(\CA_\qct(\FX))$, it
is a pair $(\CE,\alpha)$ where $\CE \in \D(\CA_\qct(\FX))$ and $\alpha: \CE
\to \CO'_\FX$ is such that $\id_\CE \otimes^{\LL}_{\CO_\FX} \alpha$ and
$\alpha
\otimes^{\LL}_{\CO_\FX} \id_\CE$ are equal isomorphisms from 
$\CE \otimes^{\LL}_{\CO_\FX} \CE$ to $\CE$.

\begin{cor}
There is a bijective correspondence between $\otimes$-compatible localizations
and idempotent pairs in $\D(\CA_\qct(\FX))$.
\end{cor}

\begin{proof} A $\otimes$-compatible localization associated to
the stable for specialization subset $Z$ gives an idempotent pair
$(\R{}\varGamma_Z(\CO'_\FX), t)$ with $t:\R{}\varGamma_Z(\CO'_\FX) \to
\CO'_\FX$ the canonical map. The condition that $\id_{}
\otimes^{\LL}_{\CO_\FX} t$ and $t \otimes^{\LL}_{\CO_\FX} \id_{}$ are equal
isomorphisms is simply the fact that $\R{}\varGamma_Z$ is an acyclization
functor associated to a $\otimes$-compatible localization.

Given an idempotent pair $(\CE,\alpha)$, define the endofunctor $\gmm$
by $\gmm(\CF) := \CF \otimes^{\LL}_{\CO_\FX} \CE$ and analogously for
morphisms. The idempotence of $\gmm$ follows from the the condition of
idempotent pair, which also ensures that it is $\otimes$-compatible. These
constructions are mutually inverse because if $Z \subset \FX$ is the stable
for specialization subset associated to $\gmm$, then 
$\R{}\varGamma_Z(\CO'_\FX) = \gmm(\CO'_\FX) 
= \CO'_\FX  \otimes^{\LL}_{\CO_\FX} \CE = \CE$.
\end{proof}  

For a complex $\CF \in \D(\CA_\qct(\FX))$ we define its \emph{homological
support} as the union of the supports of its homologies, \ie 
$\suph{\CF} := \bigcup_{i \in \ZZ} \supp{\CH^i(\CF)}$. Note that $\suph{\CF}$
is always a subset of $\FX$ stable for specialization. In fact it can be
characterized in terms of cohomology with supports, as the following result
shows.

\begin{thm} \label{suph} Let $Z \subset \FX$ be a stable for specialization
subset, for $\CF \in \D(\CA_\qct(\FX))$, we have the following equivalent
conditions:
    \begin{enumerate}
        \item $\R{}\varGamma_Z \CF  \simeq \CF$.
        \item $\CF \in \CL_Z$.
        \item $\suph{\CF} \subset Z$.
    \end{enumerate}
\end{thm}

\begin{proof} The equivalence (\textit{i}) $\dimp$ (\textit{ii}) follows from
the fact that $\CL_Z$ is a localizing subcategory with associated Bousfield
acyclization $\R{}\varGamma_Z$ as is explained in the example on page
\pageref{sfs}. The implication (\textit{i}) $\imp$ (\textit{iii}) is clear
because $\suph{\R{}\varGamma_Z \CF} \subset Z$, as $\R{}\varGamma_Z \CF$ is
computed by a complex formed by sheaves already supported in $Z$.

Let us show then that (\textit{iii}) $\imp$ (\textit{ii}). By Corollary
\ref{crithaz} it is enough to check that $\CK(x) \otimes^{\LL}_{\CO_\FX} \CF
\in \CL_Z$, for all $x \in \FX$. If $x \in Z$ then 
$\CK(x) \otimes^{\LL}_{\CO_\FX} \CF \in \CL_x \subset \CL_Z$. For $x \notin
Z$, $\FX_x \cap Z = \emptyset$ because $Z$ is stable for specialization. Let
us consider the chain of isomorphisms:
\[ \CK(x) \otimes^{\LL}_{\CO_\FX} \CF \simeq
\R{}i_{x*} i^*_x \CK(x) \otimes^{\LL}_{\CO_\FX} \CF
\overset{(\ref{compL})}{\simeq}
\CK(x) \otimes^{\LL}_{\CO_\FX} \R{}i_{x*} i^*_x \CF.
\]
Note that $\R{}i_{x*} i^*_x \CF = 0$ because 
$\suph{\CF} \subset Z \subset \FX \setminus \FX_x$, therefore we conclude that
$\CK(x) \otimes^{\LL}_{\CO_\FX} \CF  = 0$.
\end{proof}

This last result allows us to compare our classification of
$\otimes$-compatible localizations with Thomason's localization. It
says (\cite[Theorem 3.15]{th}) that there is a bijection between the set
of subsets stable for specialization of a quasi-compact quasi-separated
scheme $X$ and the set of thick triangulated $\otimes$-subcategories of
$\D(\CA_\qc(X))_\pf$. We recall that a triangulated subcategory $\CB
\subset \D(\CA_\qc(X))_\pf$ is called \emph{thick} if it is stable for direct
summands and is called by Thomason a $\otimes$-subcategories if it is a
$\otimes$-ideal, \ie the same condition that we use to define rigid
localizing subcategories. If $X$ is noetherian and separated we are able to
compare this classification with ours, which is expressed in Corollary
\ref{clastens}. We have the following:

\begin{prop} Let $X$ be a noetherian separated scheme. There is a bijection
between the set of $\otimes$-compatible localizing subcategories of 
$\D(\CA_\qc(X))$ and the set of thick triangulated
$\otimes$-subcategories of $\D(\CA_\qc(X))_\pf$. This bijection is
compatible with the classification of both sets in terms of stable
for specialization subsets of $X$.
\end{prop}
 
\begin{proof}Denote by $\loc_\otimes{(\D(\CA_\qc(X)))}$ the set of
$\otimes$-compatible localizing subcategories of
$\D(\CA_\qct(\FX))$ and by $\thik_\otimes{(\D(\CA_\qc(X))_\pf)}$ the set of
thick triangulated $\otimes$-subcategories of
$\D(\CA_\qc(X))_\pf$. Let us define a couple of maps:
\[
\loc_\otimes{(\D(\CA_\qc(X)))}
\underset{g}{\overset{f}{\rightleftarrows}}
\thik_\otimes{(\D(\CA_\qc(X))_\pf)}.
\] and check that they are mutual inverses. For a $\otimes$-compatible
localizing subcategory $\CL$ we define $f(\CL) := \CL \cap
\D(\CA_\qc(X))_\pf$ which is clearly a thick triangulated
$\otimes$-subcategory. For such a subcategory $\CB$ we define $g(\CB)$ as the
smallest localizing subcategory $\CL(\CB)$ of $\D(\CA_\qc(X))$ that contains
$\CB$. Let us show that $\CL(\CB)$ is $\otimes$-compatible. For $\CN \in
\D(\CA_\qc(X))_\pf$, define $\CL_0 \:= \{\CM \in \CL(\CB)  \,/\, 
\CM \otimes^{\LL}_{\CO_\FX} \CN \in \CL(\CB) \}$. Note that $\CL_0$ is a
localizing subcategory of $\D(\CA_\qc(X))$ and that
$\CB \subset \CL_0 \subset \CL(\CB)$, so $\CL_0 = \CL(\CB)$. Therefore,
$\CL' := \{\CN \in \D(\CA_\qc(X))  \,/\, 
\CM \otimes^{\LL}_{\CO_\FX} \CN \in \CL(\CB), \forall \CM \in \CL(\CB) \}$ is
a localizing subcategory of $\D(\CA_\qc(X))$ that contains
$\D(\CA_\qc(X))_\pf$. Applying \cite[Proposition 2.5]{Ngd}, we conclude that
$\CL' = \D(\CA_\qc(X))$, therefore $\CL(\CB)$ is rigid. The coproduct of
$\CL(\CB)$-local objects is again $\CL(\CB)$-local because $\CL(\CB)$ is
generated by perfect complexes. Then, $\CL(\CB)$ is $\otimes$-compatible by
Theorem \ref{tenscomp}.

First, let us see that $f(g(\CB)) = \CB$.  By the cited Thomason's result
there is a stable for specialization subset $Z$ of $X$ such that $\CB$ is the
class of all perfect complex with homological support contained in $Z$. It
follows that the smallest localizing subcategory that contains $\CB$,
$\CL(\CB)$, is contained in $\CL_Z$ because all of is complexes are supported
in $Z$ by Theorem \ref{suph}. Now $\CL(\CB)$ is $\otimes$-compatible, so
there is a stable for specialization subset $Z' \subset Z$ of $X$ such that 
$\CL(\CB) = \CL_{Z'}$. But $Z'$ has to agree with $Z$, otherwise by
\cite[Lemma 3.4]{th} we could find a perfect complex in $\CB$ with homological
support outside $Z'$, a contradiction. So, necessarily $\CL(\CB) = \CL_Z$ and
$\CL(\CB) \cap \D(\CA_\qc(X))_\pf = \CB$. 

Take now a $\otimes$-compatible localizing subcategory $\CL \subset
\D(\CA_\qc(X))$. By Corollary \ref{clastens}, there is a subset $Z \subset
\FX$ stable for specialization such that $\CL = \CL_Z$ which means that the
objects in $\CB := \CL \cap \D(\CA_\qc(X))_\pf$ are perfect complexes whose
homological support is contained in $Z$. The localizing subcategory $\CL' :=
g(f(\CL))$ is the smallest one that contains the objects of $\CB$, so $\CL'
\subset \CL$. The localizing subcategory $\CL'$ is $\otimes$-compatible, then
there is a stable for specialization subset $Z' \subset Z$ of $X$ such that 
$\CL' = \CL_{Z'}$. But observe that $Z'$ has to agree with $Z$ arguing as
before with the perfect complexes in the thick $\otimes$-subcategory $f(\CL)$.
\end{proof}

\begin{cor} In the previous situation, a $\otimes$-compatible localizing
subcategory of $\D(\CA_\qc(X))$ is generated by perfect complexes.
\end{cor} 
 
\begin{cosa}
Let $\CL$ be a rigid localizing subcategory of $\D(\CA_\qct(\FX))$ and $\CF,
\CG \in \D(\CA_\qct(\FX))$. The morphism
$\dhom^\cdot_{\FX}(\CF, \CG) \to 
\dhom^\cdot_{\FX}(\CF, \ell \CG)$ induced by $ \CG \to \ell \CG$
factors through $\ell\dhom^\cdot_{\FX}(\CF, \CG)$ by Proposition
\ref{righom}. So, it gives a natural morphism
\[ q: \ell\dhom^\cdot_{\FX}(\CF, \CG) \lto
      \dhom^\cdot_{\FX}(\CF, \ell \CG). 
\] Let us denote by
\[ h: \gmm\!\dhom^\cdot_{\FX}(\CF, \CG)  \lto
      \dhom^\cdot_{\FX}(\CF, \gmm \CG) 
\] the morphism such that the diagram
\begin{footnotesize}
\[\begin{CD}
\dhom^\cdot_{\FX}(\CF, \gmm \CG)  
             @>>> \dhom^\cdot_{\FX}(\CF, \CG)
             @>>> \dhom^\cdot_{\FX}(\CF, \ell \CG) @>+>>     \\ 
         @AAhA                @|             @AAqA             \\ 
\gmm\!\dhom^\cdot_{\FX}(\CF, \CG) 
             @>>> \dhom^\cdot_{\FX}(\CF, \CG)
             @>>> \ell\dhom^\cdot_{\FX}(\CF, \CG) @>+>>
\end{CD}
\]\end{footnotesize} is a morphism of distinguished triangles. Again, $h$
and $q$ determine each other.

With the notation of the previous remark, we say that the
localization $\ell$ is $\shom$\emph{-compatible} (or that $\CL$ is
$\shom$-compatible or that $\gmm$ is $\shom$-compatible) if the
canonical morphism $q$, or equivalently $h$, is an isomorphism.
\end{cosa} 

\begin{cosa} Let $\CL_Z$ be a $\otimes$-compatible localizing
subcategory of $\D(\CA_\qct(\FX))$ whose associated (stable for
specialization) subset is $Z \subset \FX$. Let us apply the functor 
$\dhom^\cdot_{\FX}(-, \CF)$, where $\CF \in \D(\CA_\qct(\FX))$, to the
canonical triangle
\[ \gmm_{\!Z} \CO'_\FX \lto \CO'_\FX \lto \ell_Z \CO'_\FX \overset{+}{\lto}
\] associated to $\CL_Z$.
We have added the associated subsets as subindices for clarity. We obtain:
\begin{equation}\label{trihom}
\dhom^\cdot_{\FX}(\ell_Z \CO'_\FX, \CF) \lto \CF
    \lto \dhom^\cdot_{\FX}(\gmm_{\!Z} \CO'_\FX, \CF)
\overset{+}{\lto}.
\end{equation} \end{cosa} 

\begin{prop} the canonical natural transformations:
\[\id \to \dhom^\cdot_{\FX}(\gmm_{\!Z} \CO'_\FX, -)
\text{ \quad and \quad }
\dhom^\cdot_{\FX}(\ell_Z \CO'_\FX, -) \to \id,
\] 
correspond to a $\shom$-compatible localization and its
corresponding acyclization in $\D(\CA_\qct(\FX))$, respectively. Its
associated subset of $\FX$ is $\FX \setminus Z$.
\end{prop} 

\begin{proof} Note that (\ref{trihom}) is a Bousfield localization triangle
because $\CL_Z$ is $\otimes$-compatible.  The associated localizing
subcategory
\[ \CL  = \{ \CM \in \D(\CA_\qct(\FX)) \,/\,  
   \dhom^\cdot_{\FX}(\gmm_{\!Z} \CO'_\FX, \CM) = 0\}
\]
satisfies that ${^\perp{(\CL^\perp)}} = \CL$ (\cite[Proposition 1.6]{AJS}).
Furthermore, the canonical isomorphisms 
\begin{equation*}
   \begin{split}
      \dhom^\cdot_{\FX}(\CF, \dhom^\cdot_{\FX}(\gmm_{\!Z} \CO'_\FX, \CG)) 
   &\cong \dhom^\cdot_{\FX}(\CF \otimes^{\LL}_{\CO_\FX} \gmm_{\!Z} \CO'_\FX,
\CG)
           \\
   &\cong \dhom^\cdot_{\FX}(\gmm_{\!Z} \CO'_\FX, \dhom^\cdot_{\FX}(\CF, \CG))
   \end{split} 
\end{equation*}
show that $\CL$ is rigid (Proposition \ref{righom}) and $\shom$-compatible.

Let us check that $\CL = \CL_{\FX \setminus Z}$. Let $z \in \FX$,
we will consider two possibilities depending on the point being or not in $Z$.
First, if $z \in \FX \setminus Z$, it follows that
$\CK(z) \in \CL_Z^{\perp}$ by Lemma \ref{locKhom} and therefore we have
that
\[\dhom^\cdot_{\FX}(\gmm_{\!Z} \CO'_\FX, \CK(z)) \osi 
  \dhom^\cdot_{\FX}(\gmm_{\!Z} \CO'_\FX, \gmm_{\!Z} \CK(z)) =0,\] 
then
\[\dhom^\cdot_{\FX}(\ell_Z \CO'_\FX, \CK(z)) \iso \CK(z).\]

For $z \in Z$ we will show that $\dhom^\cdot_{\FX}(\ell_Z \CO'_\FX, \CK(z)) =
0$. By Proposition \ref{smalloc} it is enough to prove that
\[ \Hom_{\D(\FX)}(\CK(y), \dhom^\cdot_{\FX}(\ell_Z \CO'_\FX, \CK(z))) = 0,
 \quad \forall y \in \FX,
\]
equivalently that
\[ \Hom_{\D(\FX)}(\CK(y) \otimes^{\LL}_{\CO_\FX} \ell_Z \CO'_\FX, \CK(z)) = 0,
 \quad \forall y \in \FX.
\]
The localization functor $\ell_Z$ is $\otimes$-compatible so 
$\CK(y) \otimes^{\LL}_{\CO_\FX} \ell_Z \CO'_\FX \cong \ell_Z \CK(y)$ will be
zero if $y \in Z$. On the other hand, if $y \in \FX \setminus Z$ we conclude
because $\CK(y) \otimes^{\LL}_{\CO_\FX} \ell_Z \CO'_\FX \in \CL_y$ and
$\CK(z) \in \CL_y^{\perp}$ (Lemma \ref{locKhom}).
\end{proof} 

\begin{cosa} Note that the following adjunction is completely formal
\[\dhom^\cdot_{\FX}(\gmm_{\!Z} \CF, \CG) \liso
\dhom^\cdot_{\FX}(\CF, \ell_{\FX \setminus Z}\CG).
\]
Indeed, it is the composition of the following natural isomorphisms
\begin{equation*}
   \begin{split}
    \dhom^\cdot_{\FX}(\gmm_{\!Z} \CF, \CG)
          & \cong 
              \dhom^\cdot_{\FX}(\gmm_{\!Z} \CO'_\FX
                            \otimes^{\LL}_{\CO_\FX} \CF, \CG) \\
          & \cong
              \dhom^\cdot_{\FX}(\CF,
               \dhom^\cdot_{\FX}(\gmm_{\!Z} \CO'_\FX, \CG)) \\
          & \cong \dhom^\cdot_{\FX}(\CF,\ell_{\FX \setminus Z}\CG).
   \end{split}
\end{equation*}
\end{cosa}

\begin{ex} Let now $Z$ be a closed subset of $X$ which we assume it is an
ordinary (noetherian separated) scheme. and $\LL{}\varLambda_Z:
\D(\CA_\qc(X)) \to \D(\CA(X))$ the left-derived functor of the
completion along the closed subset $Z$ (which exist because it can be
computed using quasi-coherent flat resolutions, as proved in \cite{AJL1}). In
\emph{loc.~cit.} it is also shown there is a natural isomorphism
\[\dhom^\cdot_{X}(\R{}\varGamma_Z \CO_X, \CG) \liso 
\R{}Q\LL{}\varLambda_Z (\CG).
\] This result together with the previous adjunction is often referred to as
\emph{Greenlees-May duality} because it generalizes a result from \cite{gm}
in the affine case.
\end{ex}

\begin{cosa} In general, if $Z \in \FX$ is a stable for specialization
subset of $\FX$ we will define for every $\CG \in \D(\CA_\qct(\FX))$:
\[\boldsymbol{\Lambda}_Z(\CG) := 
  \dhom^\cdot_{\FX}(\R{}\varGamma_Z \CO'_\FX, \CG)).
\]
\end{cosa}

\begin{thm} \label{homcomp} For a rigid localizing subcategory $\CL \subset
\D(\CA_\qct(\FX))$, the following are equivalent:
    \begin{enumerate}
        \item \label{hc1} The localization associated to $\CL$ is  
                          $\shom$-compatible.
        \item \label{hc1'} For every $\CN \in \CL$ and 
                          $\CF \in \D(\CA_\qct(\FX))$ we have that
                          $\dhom^\cdot_{\FX}(\CF, \CN) \in \CL$.
        \item \label{hc4} The set $Y := \psi(\CL)$ is generically
                          stable\footnote{\ie an arbitrary intersection of
                          open subsets.} and its associated localization
                          functor is $\boldsymbol{\Lambda}_Z$ with 
                          $Z = \FX \setminus Y$.
    \end{enumerate}

\end{thm}

\begin{proof} 
Let us see first that (\textit{\ref{hc1'}}) $\imp$ (\textit{\ref{hc4}}). 
Let $z \in Y$ and $x \in \FX$ such that $z \in \overline{\{x\}}$. With the
notation of Corollary \ref{capiny}, if $x \notin Y$ then by Lemma
\ref{locKhom2}, $\dhom^\cdot_{\FX}(\CE(x), \CE(z)) \in \CL^{\perp}$.
By (\textit{\ref{hc1'}}), $\dhom^\cdot_{\FX}(\CE(x), \CE(z))$ belongs to
$\CL$ because $\CE(z) \in \CL$. Therefore $\dhom^\cdot_{\FX}(\CE(x), \CE(z))
= 0$ and we have
\[ \Hom_{\D(\FX)}(\CE(x), \CE(z)) \cong
\Hom_{\D(\FX)}(\CO'_\FX,\dhom^\cdot_{\FX}(\CE(x), \CE(z))) = 0
\] a contradiction. Necessarily, the set $Z = \FX \setminus Y$ is stable for
specialization and $\ell_Y = \boldsymbol{\Lambda}_Z$.

The implication (\textit{\ref{hc4}}) $\imp$ (\textit{\ref{hc1}}) follows from
the previous remarks and the bijective correspondence established in
Theorem \ref{ThclasB}.

To finish, (\textit{\ref{hc1}}) $\imp$ (\textit{\ref{hc1'}}) is
straightforward because for every $\CN \in \CL$, we have that
\[ \dhom^\cdot_{\FX}(\CF, \CN) = \dhom^\cdot_{\FX}(\CF, \gmm \CN)
      \overset{\text{(\textit{\ref{hc1}})}}{\cong} 
      \gmm \dhom^\cdot_{\FX}(\CF, \CN) \in \CL
\]
\end{proof}

\begin{cor} The functor $\gmm$ associated to a $\shom$-compatible
localization in $\D(\CA_\qct(\FX))$ commutes with \emph{products}, in
particular, the corresponding localizing class $\CL$ is closed for products.
\end{cor}

\begin{proof} It is an immediate consequence of Theorem \ref{homcomp} 
(\textit{\ref{hc4}}) and that every complex in $\K(\CA_\qct(\FX))$ admits a
K-flat resolution by \cite[Proposition 2.1.3]{fgm} and a K-injective
resolution.
\end{proof}  

\begin{cor} For a noetherian separated formal scheme $\FX$ there is a
bijection between the class of $\shom$-compatible localizations of
$\D(\CA_\qct(\FX))$ and the set of generically stable subsets of $\FX$.
\end{cor}

\end{document}